\documentclass[10pt]{article}

\usepackage{amsmath,amssymb,fullpage}

\usepackage{textcomp,marvosym}

\usepackage{cite}

\usepackage{hyperref}

\usepackage[table]{xcolor}

\usepackage{array}

\usepackage{here}
\usepackage{graphicx}

  \usepackage{zref-xr}
  \zxrsetup{tozreflabel=false, toltxlabel=true, verbose}
\zexternaldocument*{./StochasticPI_PLOS_SI}

\def\Xz{\boldsymbol{X}}

\providecommand\phib{\boldsymbol{\emptyset}}
\providecommand\X[1]{\boldsymbol{X_{#1}}}
\providecommand\Z[1]{\boldsymbol{Z_{#1}}}

\def\E{\mathbb{E}}
\def\V{ \textnormal{Var} }
\def\C{ \textnormal{Cov} }

\def\llongrightarrow{\relbar\joinrel\relbar\joinrel\relbar\joinrel\rightarrow}

\providecommand{\rarrow}[1]{\stackrel{#1}{\llongrightarrow}}

\DeclareMathOperator*{\col}{col}

\title{Variance reduction for antithetic integral control of stochastic reaction networks}
\author{Corentin Briat\thanks{These two coauthors contributed equally}, Ankit Gupta\footnotemark[1], Mustafa Khammash\thanks{Corresponding author: mustafa.khammash@bsse.ethz.ch}\\ \ \\Department of Biosystems Science and Engineering,\\ETH-Z\"urich, Switzerland}
\date{}

\begin{document}

\maketitle
\section*{Abstract}

The antithetic integral feedback motif recently introduced in \cite{Briat:15e} is known to ensure robust perfect adaptation for the mean dynamics of a given molecular species involved in a complex stochastic biomolecular reaction network. However, it was observed that it also leads to a higher variance in the controlled network than that obtained when using a constitutive (i.e. open-loop) control strategy. This was interpreted as the cost of the adaptation property and may be viewed as a performance deterioration for the overall controlled network. To decrease this variance and improve the performance, we propose to combine the antithetic integral feedback motif with a negative feedback strategy. Both theoretical and numerical results are obtained. The theoretical ones are based on a tailored moment closure method allowing one to obtain approximate expressions for the stationary variance for the controlled network and predict that the variance can indeed be decreased by increasing the strength of the negative feedback. Numerical results verify the accuracy of this approximation and show that the controlled species variance can indeed be decreased, sometimes below its constitutive level. Three molecular networks are considered in order to verify the wide applicability of two types of negative feedback strategies. The main conclusion is that there is a trade-off between the speed of the settling-time of the mean trajectories and the stationary variance of the controlled species; i.e. smaller variance is associated with larger settling-time.

\section*{Author summary}

Homeostasis, the ability of living organisms to regulate their internal state, is of fundamental importance for their adaptation to environmental changes and their survival. This is the reason why complex regulatory genetic networks evolved and allowed for the emergence of more and more complex organisms. Recently, the theoretical study of those regulatory networks using ideas and concepts from control theory and the design of novel ones have gained a lot of attention. Synthetic regulatory circuits are seen as elementary building blocks for the design of complex networks that need to incorporate some regulating elements to be fully functional. This is for instance the case in metabolic engineering where the production of biomolecules, such as drugs or biofuels, need to be optimized and tightly regulated. A particular circuit, the so-called antithetic integral controller, is now known to ensure homeostasis even when regulatory circuits are subject to randomness. However, it is also known that this circuit increases variability in the network. The effects of a correcting negative feedback loop on the variance are discussed here and it is shown that variability can be reduced this way. Notably, we show that there is a tradeoff between speed of the network and variability.

\section*{Introduction}

The design and implementation of artificial in-vivo biomolecular controllers have become very popular \cite{Briat:15e,Briat:16a,Qian:17,Cuba:17,Annunziata:17,Lillacci:17} because of their potential applications for the tight and robust control of gene expression \cite{Briat:15e}, the optimization of metabolic networks for the efficient production of biomolecules \cite{Venayak:15,Cress:15} , or the development of new treatments for certain genetic diseases \cite{Ye:14}. Indeed, many of the instances of those problems can be interpreted from an homeostatic point of view in the sense that they may all be solved by achieving or restoring homeostasis in the corresponding genetic network using synthetic regulatory circuits \cite{Ye:14,Venayak:15,Cress:15,Briat:15e,Schukur:16}. In this regard, those problems essentially reduce to the design and the implementation of robust and reliable regulatory circuits that can optimize an inefficient network or correct a malfunctioning one -- an observation which strongly suggests that ideas from control theory and control engineering \cite{Albertos:10} could be adapted to biochemical control problems \cite{DelVecchio:15,Harris:15,Briat:15e}. A cornerstone in control theory and engineering is the so-called \emph{integral controller} that can ensure precise constant set-point regulation for a regulated variable in a given system. Such mechanism, where the action onto the controlled system is depending on the integral of the deviation of the regulated variable from the desired set-point, is to be contrasted with the so-called \emph{proportional controller} where the system is simply actuated proportionally to the deviation of the regulated variable from the desired set-point. Unlike integral control, the latter one is unable to achieve robust constant set-point regulation for the controlled variable and to reject constant disturbances. In other words, integral control has the capacity of ensuring perfect adaptation for the regulated variable. The downside, however, is that it has a destabilizing effect on the dynamics (emergence of oscillations or even diverging trajectories) of the overall controlled system which can be then compensated by adjoining a proportional action, thus giving rise to the so-called Proportional-Integral (PI) controller \cite{Astrom:95}.

Based on the strength of these facts, an integral controller referred to as the \emph{antithetic integral controller} was proposed in \cite{Briat:15e} for the control of the mean level of some molecular species of interest in a given biochemical reaction network. This controller is implementable in terms of elementary biochemical reactions with mass-action kinetics, making it practically implementable in-vivo using, for instance, sigma- and anti-sigma-factors \cite{Lillacci:17}. This controller theoretically works in both the deterministic and the stochastic settings. In the latter setting, it was notably shown that, under some reasonable conditions, the ergodicity properties of the controlled network are independent of the parameters of the antithetic integral controller -- a surprising key property that has no counterpart in the deterministic setting and that dramatically simplifies its implementation. A drawback, however, is the increase of the stationary variance of the regulated species compared to the constitutive variance that would be obtained by using a static open-loop strategy, even though the latter one would be unable to ensure regulation and perfect adaptation for the mean level of the regulated species. This phenomenon is seemingly analogous to the destabilizing behavior of the deterministic integral controller mentioned in the previous paragraph. This variance increase can hence be interpreted as the price to pay for perfect adaptation at the mean species level.

The goal of this paper is to investigate the effect of adding a negative feedback to the antithetic integral motif in a way akin, yet different, to deterministic PI controllers. As discussed above, adding a proportional action in the deterministic setting compensates for the destabilizing effect of the integrator. Comparatively, it may seem reasonable to think that, in the stochastic setting, a proportional action could have an analogous effect and would result in a decreased variance  for the controlled variable (this is, for instance, what happens when considering certain linear systems driven by white noise). In fact, it has been shown that negative feedback at a transcriptional level in a gene expression network leads to a variance reduction in the protein levels; see e.g. \cite{Becskei:00,Thattai:01,Paulsson:04,Kaern:05} and the references therein.

Two types of negative feedback are considered in the present paper: the first one consists of an ON/OFF proportional action whereas the second one is governed by a repressing Hill function. First we theoretically prove using a tailored moment closure method that, in a gene expression network controlled with an antithetic integral controller, the stationary variance in the protein copy number can be decreased by the use of a negative feedback. In this specific case, the steady-state variance is decreasing monotonically as a function of the strength of the negative feedback. An immediate consequence is that it is theoretically possible to set the steady-state variance to a level that lies below the constitutive steady-state variance, which is the value of the steady-state variance that would have been obtained using a constitutive (i.e. open-loop) control strategy. The theoretical prediction will also be observed by exact numerical predictions using Gillespie's algorithm (Stochastic Simulation Algorithm - SSA \cite{Gillespie:76}). A caveat, however, is that setting the gain of the negative feedback very high will likely result in a very low steady-state variance but may also result in a regulation error for the mean of the controlled species and in a loss of ergodicity for the overall controlled network. In this regard, reducing the steady-state variance below its constitutive level may not always be physically possible. Finally, it is also emphasized that a low stationary variance for the controlled species is often associated with higher settling-time for the controlled species. Hence, there is a tradeoff between variability and fast dynamics/small settling-time. The two negative feedback actions also exhibit quite different behaviors. Indeed, while the ON/OFF proportional feedback seems to be efficient at reducing the stationary variance through an increase of its gain, the dynamics of the mean gets first improved by reducing the settling-time but then gets dramatically deteriorated by the appearance of a fast initial transient phase followed by a very slow final one resulting then in a high settling-time. On the other hand, the Hill controller leads to very homogeneous mean dynamics for different feedback strength but the steady-state variance is also much less sensitive and does not vary dramatically. It is then argued that those differences may find an explanation by reasoning in a deterministic point of view. The ON/OFF controller (an error-feedback) introduces a stable zero in the dynamics of the closed-loop network which is small in magnitude when the gain of the negative feedback is high. When this is the case,  the zero is close to the origin and the closed-loop dynamics will almost contain a derivative action, whence the fast initial transient phase. On the other hand, the Hill negative feedback (an output-feedback) does not introduce such a zero in the closed-loop dynamics, which may explain the homogeneity of the mean trajectories. Another possible reason is that the effective proportional gain (which will be denoted by $\beta$) is much less sensitive to changes in the feedback strength than the ON/OFF controller.

Approximate equations for the stationary variance are then obtained in the general unimolecular network case. The obtained expressions shed some light on an interesting connection between the covariances of the molecular species involved in the stochastic reaction network and the stability of a deterministic linear system controlled with a standard PI controller, thereby unveiling an unexpected, yet coherent, bridge between the stochastic and deterministic settings. Applying this more general framework to the a gene expression network with protein maturation allows one to reveal that the steady-state variance may not be necessarily a monotonically decreasing function of the negative feedback strength. In spite of this, the same conclusions as in the gene expression network hold: the variance can sometimes be decreased below its constitutive level but this may also be accompanied with a loss of ergodicity. The same qualitative conclusions for the transient of the mean dynamics and the properties of the controller also hold in this case.

Even though the proposed theory only applies to unimolecular networks, stochastic simulations are performed for a gene expression network with protein dimerization; a bimolecular network. Once again, the same conclusions as in for previous networks hold with the difference that the constitutive variance level is unknown in this case due to openness of the moment equations. These results tend to suggest that negative feedback seems to operate in the same way in bimolecular networks as in unimolecular networks.


\subsection*{Reaction networks}

Let us consider a stochastic reaction network $(\Xz,\mathcal{R})$ involving $d$ molecular species $\X{1},\ldots,\X{d}$ interacting through $K$ reaction channels $\mathcal{R}_1,\ldots,\mathcal{R}_K$ defined as
\begin{equation}
 \mathcal{R}_k:\ \sum_{i=1}^d\zeta_{k,i}^l\X{i}\rarrow{\rho_k}  \sum_{i=1}^d\zeta_{k,i}^r\X{i},\ k=1,\ldots,K
\end{equation}
where $\rho_k\in\mathbb{R}_{>0}$ is the reaction rate parameter and $\zeta_k^r=\col(\zeta_{k,1}^r,\ldots,\zeta_{k,d}^r)$, $\zeta_k^l=\col(\zeta_{k,1}^l,\ldots,\zeta_{k,d}^l)$ are the left/right stoichiometric vectors of the reaction $\mathcal{R}_k$. The corresponding stoichiometric vector is hence given by $\zeta_k:=\zeta_k^r-\zeta_k^l\in\mathbb{Z}^d$ indicating that when this reaction fires, the state jumps from $x$ to $x+\zeta_k$. The stoichiometric matrix $S\in\mathbb{Z}^{d\times K}$ of this reaction network is defined as $S:=\begin{bmatrix}
  \zeta_1\ \cdots\ \zeta_K
\end{bmatrix}$. When the kinetics is mass-action, the propensity function $\lambda_k$ of the reaction $\mathcal{R}_k$ is given by $\textstyle\lambda_k(x)=\rho_k\prod_{i=1}^d\frac{x_i!}{(x_i-\zeta_{n,i}^l)!}$. Under the well-mixed assumption, the above network can be described by a continuous-time Markov process $(X_1(t),\ldots,X_d(t))_{t\ge0}$ with the $d$-dimensional nonnegative lattice $\mathbb{Z}_{\ge0}^d$ as state-space; see e.g. \cite{Anderson:11}.

\subsection*{The regulation/perfect adaptation problems and antithetic integral control}

Let us consider here a stochastic reaction network  $(\Xz,\mathcal{R})$. The regulation problem consists of finding another reaction network (i.e. a set of additional species and additional reactions) interacting with $(\Xz,\mathcal{R})$ in a way that makes the interconnection well-behaved (i.e. ergodic) and such that the mean of some molecular species $\X{\ell}$ for some given $\ell\in\{1,\ldots,d\}$ converges to a desired set-point (given here by $\mu/\theta$ for some $\mu,\theta>0$) in a robust way; i.e. irrespective of the values of the parameters of the network $(\Xz,\mathcal{R})$.

It was shown in \cite{Briat:15e} that, under some assumptions on the network  $(\Xz,\mathcal{R})$, the antithetic integral controller defined as
\begin{equation}\label{eq:AIC}
  \underbrace{\phib\rarrow{\mu}\Z{1}}_{\mbox{reference}},\ \underbrace{\phib\rarrow{\theta X_\ell}\Z{2}}_{\mbox{measurement}},\ \underbrace{\Z{1}+\Z{2}\rarrow{\eta}\phib}_{\mbox{comparison}},\ \underbrace{\phib\rarrow{kZ_1}\X{1}}_{\mbox{actuation}},
\end{equation}
solves the above regulation problem. This regulatory network consists of two additional species $\Z{1}$ and $\Z{2}$, and four additional reactions. The species $\Z{1}$ is referred to as the actuating species as it is the species that governs the rate of the actuation reaction which produces the actuated species $\X{1}$ at a rate proportional to $Z_1$. The species $\Z{2}$ is the sensing species as it is produced at a rate proportional to the controlled species $\X{\ell}$ through the measurement reaction. The first reaction is the reference reaction as it encodes part of the set-point $\mu/\theta$ whereas the third reaction is the comparison reaction that compares the population of the controller species and annihilates them accordingly, thereby closing negatively the loop while, and at the same time correlating the populations of the controller species. The comparison (or titration) reaction is the crucial element of the above controller network and, to realize such a reaction, one needs to rely on intrinsic strongly binding properties of certain molecules such as sigma- and anti-sigma-factors \cite{Briat:15e} or small RNAs and RNAs \cite{Qian:17,Levine:07,Yoo:13}.

\subsection*{Variance amplification in antithetic integral control}

We discussed above about the convergence properties of the mean level of the controlled species $\X{\ell}$ when network $(\Xz,\mathcal{R})$ is controlled with the antithetic integral controller \eqref{eq:AIC}. However, it was remarked in \cite{Briat:15e} that while the mean of  $\X{\ell}$ converges to the desired steady-state, the stationary variance of the controlled species could be much larger than its constitutive value that would be obtained by simply considering a naive constitutive production of the species $\X{1}$ that would lead to the same mean steady-state value $\mu/\theta$. This was interpreted as the price to pay for having the perfect adaptation property for the controlled species. To illustrate this phenomenon, let us consider the following gene expression network:
\begin{equation}\label{eq:gene_expression}
  \phib\rarrow{k_r}\X{1},\ \X{1}\rarrow{k_p}\X{1}+\X{2}, \X{1}\rarrow{\gamma_r}\phib, \X{2}\rarrow{\gamma_p}\phib
\end{equation}
where $\X{1}$ denotes mRNA and $\X{2}$ denotes protein. The objective here is to control the mean level of the protein by acting at a transcriptional level using the antithetic controller \eqref{eq:AIC}; hence, we set $k_r=kZ_1$. Using a tailored moment closure method, it is proved in the SI that the stationary variance $\V_\pi^{I}(X_2)$ for the protein copy number is approximately given by the following expression
\begin{equation}\label{eq:variance_ge}
  \V_\pi^{I}(X_2)\approx\dfrac{\mu}{\theta}\left(\dfrac{1+\dfrac{k_p}{\gamma_r+\gamma_p}+\dfrac{kk_p}{\gamma_r\gamma_p}}{1-\dfrac{k\theta k_p}{\gamma_r\gamma_p(\gamma_r+\gamma_p)}}\right),\ k>0,\ k/\eta\ll1.
\end{equation}
The rationale for the assumption $k/\eta\ll1$ is that it allows for closing the moments equation (which is open because of the presence of the comparison reaction) and obtain a closed-form solution for the stationary variance.  On the other hand, the constitutive (i.e. open-loop) stationary variance $\V_\pi^{OL}(X_2)$ for the protein copy number obtained with the constitutive strategy
\begin{equation}
  k_r=\dfrac{\mu}{\theta}\dfrac{\gamma_r\gamma_p}{k_p}
\end{equation}
is given by
\begin{equation}
  \V_\pi^{OL}(X_2)=\dfrac{\mu}{\theta}\left(1+\dfrac{k_p}{\gamma_r+\gamma_p}\right).
\end{equation}
It is immediate to see that the ratio
\begin{equation}
  \dfrac{\V_\pi^{I}(X_2)}{\V_\pi^{OL}(X_2)}\approx\dfrac{1+\dfrac{kk_p}{\gamma_r\gamma_p(\gamma_r+\gamma_p)}}{1-\dfrac{k\theta k_p(k_p+\gamma_r+\gamma_p)}{\gamma_r\gamma_p(\gamma_r+\gamma_p)}},\ k,\theta>0,\ k/\eta\ll1
\end{equation}
is greater than 1 for all $k,\theta>0$ such that the denominator is positive. Note that the above formula is not valid when $k=0$ or $\theta=0$ since this would result in an open-loop network for which set-point regulation could not be achieved. This expression is also a monotonically increasing function of the gain $k$, a fact that was numerically observed in \cite{Briat:15e}. This means that choosing $k$ very small will only result in a small increase of the stationary variance of the controlled species when using an antithetic integral feedback. However, this will very likely result in very slow dynamics for the mean of the controlled species.

Finally, it is important to stress that while this formula is obviously not valid when the denominator is nonpositive, we know from \cite{Briat:15e} that in the case of the gene expression network, the closed-loop network will be ergodic with converging first and second-order moments for all $k>0$ and all $\theta>0$ (assuming that the ratio $\mu/\theta$ is kept constant). This inconsistency stems from the fact that the proposed theoretical approach relies on a tailored moment closure approximation that will turn out to be connected to the Hurwitz stability of a certain matrix that may become unstable when the gain $k$ of the integrator is too large. This will be elaborated more in the following sections.

\subsection*{Negative feedback action}

We will consider in this paper two types of negative feedback action. The first one, referred to as the \emph{ON/OFF proportional feedback}, is essentially theoretical and cannot be exactly implemented, but it may be seen as a local approximation of some  more complex (e.g. nonlinear) repressing function. It is given by the reaction
\begin{equation}\label{eq:djksjdl}
  \phib\rarrow{F(X_\ell)}\X{1}
\end{equation}
together with the propensity function $F(X_\ell)=K_p\max\{0,\mu-\theta X_\ell\}$ where $K_p$ is the so-called feedback gain/strength. It is similar to the standard proportional feedback action used in control theory with the difference that a regularizing function, in the form of a max function, is involved in order the restrict the propensity function to nonnegative values. Note that this controller can still be employed for the in-silico control of single-cells using a stochastic controller as, in this case, we would not be restricted anymore to mass-action, Hill or Michaelis-Menten kinetics. This was notably considered in the case of in-silico population control in \cite{Briat:12c,Briat:13h,Guiver:15}.

The second type of negative feedback action, referred to as the \emph{Hill feedback}, consists of the reaction \eqref{eq:djksjdl} but involves the non-cooperative repressing Hill function $F(X_\ell)=K_p/(1+X_\ell)$ as propensity function. This type of negative feedback is more realistic as such functions have empirically been shown to arise in many biochemical, physiological and epidemiological models; see e.g. \cite{Murray:02}.

In both cases, the total rate of production of the molecular species $\X{1}$ can be expressed as the sum $kZ_1+F(X_\ell)$ which means that, at stationarity, we need to have that $\E_\pi[kZ_1+F(X_\ell)]=u^*$ where $u^*$ is equal to the value of the constitutive (i.e. deterministic) production rate for $\X{1}$ for which we would have that $\E_\pi[X_\ell]=\mu/\theta$. Noting now that for both negative feedback functions, we will necessarily have that $\E_\pi[F(X_\ell)]>0$, then this means that if the gain $K_p$ is too large, it may be possible that the mean of the controlled species do not converge to the desired set-point implying, in turn, that the overall controlled network will fail to be ergodic. This will be notably the case when $\E_\pi[F(X_\ell)]>u^*$. In particular, when $F(X_\ell)=K_p\max\{0,\mu-\theta X_\ell\}$, a very conservative sufficient condition for the closed-loop network to be ergodic is that $K_p<u^*/\mu$ whereas when $F(X_\ell)=K_p/(1+X_\ell)$, this condition becomes $K_p<u^*$. These conditions can be determined by considering the worst-case mean value of the negative feedback strategies; i.e. $K_p\mu$ and $K_p$, respectively.

\section*{Results}

\subsection*{Invariants for the antithetic integral controller}

We describe some important invariant properties of the antithetic integral controller \eqref{eq:AIC} which are independent of the parameters of the controlled network under the assumption that these invariants exist; i.e. they are finite. Those invariants are given by
\begin{align}
\label{geinv1}
\C_\pi (X_\ell,Z_1-Z_2) = \frac{\mu}{\theta},
\end{align}
\begin{align}
\label{geinv2}
\E_\pi (Z_1Z_2) = \frac{\mu}{\eta},
\end{align}
\begin{align}
\label{geinv3}
\E_\pi (Z^2_1Z_2) = \frac{\mu}{\eta} \left(  1  +  \E_\pi(Z_1) \right)
\end{align}
and
\begin{align}
\label{geinv4}
\E_\pi (Z_1Z^2_2) = \frac{\mu  +  \theta \E_\pi( X_\ell Z_2 ) }{\eta}
\end{align}
play an instrumental role in proving all the theoretical results of the paper. Interestingly, we can notice that $\C_\pi (X_\ell,Z_1-Z_2)=\E_\pi[X_\ell]$, which seems rather coincidental. From the second invariant we can observe that, if $\eta\gg\mu$, then $\E_\pi (Z_1Z_2)\approx0$, which indicates that the values taken by the random variable $Z_2(t)$ will be most of the time equal to 0. Note that it cannot be $Z_1(t)$ to be mostly taking zero values since $\Z{1}$ is the actuating species whose mean must be nonzero (assuming here that the natural production rates of the molecular species in the controlled network are small). Similarly, setting $\eta$ large enough in the third expression will lead to a similar conclusion. Note that $\E_\pi(Z_1)$ is independent of $\eta$ here and only depends on the set-point $\mu/\theta$, the integrator gain $k$ and the parameters of the network which is controlled. The last expression again leads to similar conclusions. Indeed, if $\eta$ is sufficiently large, then $\E_\pi( X_\ell Z_2 )\approx0$ and, hence, $\E_\pi (Z_1Z^2_2) \approx 0$ which implies that $Z_2(t)$ needs to be most of the time equal to 0. These properties will be at the core of the moment closure method used to obtain an approximate closed-form formula for the covariance matrix for the closed-loop network.

\subsection*{An approximate formula for the stationary variance of the controlled species}

Let us assume here that the open-loop network $(\Xz,\mathcal{R})$ is mass-action and involves, at most, unimolecular reactions. Hence, the vector of propensity functions can be written as
\begin{equation}
  \lambda(x)=Wx+w_0
\end{equation}
for some nonnegative matrix $W\in\mathbb{R}^{K\times d}$ and nonnegative vector $w_0\in\mathbb{R}^K$. It is proved in the SI that, under the assumption $k/\eta\ll1$, we can overcome the moment closure problem arising from the presence of the annihilation reaction in the antithetic controller and show that the exact stationary covariance matrix of the network given by
  \begin{equation*}
\begin{bmatrix}
\C_\pi^{PI}( X,X) & \C^{PI}_\pi( X ,Z) \\
\C_\pi^{PI}( Z, X ) & \V^{PI}_\pi(Z)
\end{bmatrix},\ Z:=Z_1-Z_2
\end{equation*}
is approximatively given by the matrix $\Sigma$ solving the Lyapunov equation
\begin{equation}\label{eq:Lyapunov}
  R\Sigma + \Sigma R^T + Q = 0
\end{equation}
where
\begin{equation*}
\begin{array}{rcl}
R&=&\begin{bmatrix}
  SW - \beta e_1 e_\ell ^T & k e_1\\
-\theta e^T_\ell   &0
\end{bmatrix},\\
Q&=&\begin{bmatrix}
S D S^T + c e_1 e^T_1 & 0 \\
0  & 2\mu
\end{bmatrix},\\
D&=&\textnormal{diag}(W\E_\pi[X]+w_0),\\
c &=& - \dfrac{1}{ e^T_\ell (SW)^{-1}e_1 } \left( \dfrac{\mu}{\theta} + e^T_\ell  (SW)^{-1}Sw_0  \right),\\
\beta&=&-\dfrac{\C^{PI}_\pi( F(X_\ell ), X_\ell)}{\C^{PI}_\pi( X_\ell , X_\ell)}.
\end{array}
\end{equation*}
Note that since the function $F$ is decreasing then the effective proportional gain, $\beta$, is always a positive constant and seems to be mostly depending on $K_p$ but does not seem to change much when $k$ varies (see e.g. Figure \ref{fig:Gene_Prop_Beta} and Figure \ref{fig:Gene_Hill_Beta} in the appendix). It can also be seen that for the Lyapunov equation to have a positive definite solution, we need that the matrix $R$ be Hurwitz stable; i.e. all its eigenvalues have negative real part. In parallel of that, it is known from the results in \cite{Briat:15e} that the closed-loop network will remain ergodic when $\beta=0$ even when the matrix $R$ is not Hurwitz stable. In this regard, the formula \eqref{eq:Lyapunov} can only be valid when the parameters $\beta$ and $k$ are such that the matrix $R$ is Hurwitz stable. When this is not the case, the formula is out its domain of validity and is meaningless. The stability of the matrix $R$ is discussed in more details in the SI.

\subsection*{Connection to deterministic proportional-integral control}

Interestingly, the matrix $R$ coincides with the closed-loop system matrix of a deterministic linear system controlled with a particular proportional-integral controller. To demonstrate this fact, let us consider the following linear system
\begin{equation}
  \begin{array}{rcl}
    \dot{x}(t)&=&SWx(t)+e_1u(t)\\
    y(t)&=&e_\ell^Tx(t)
  \end{array}
\end{equation}
where $x$ is the state of the system, $u$ is the control input and $y$ is the measured/controlled output. We propose to use the following PI controller in order to robustly steers the output to a desired set-point $\mu/\theta$
\begin{equation}
  u(t)=\dfrac{\beta}{\theta}(\mu-\theta y(t))+k\int_0^t (\mu-\theta y(s))ds
\end{equation}
where $\theta$ is the sensor gain, $\beta/\theta$ is the proportional gain and $k$ is the integral gain. The closed-loop system is given in this case by
\begin{equation}
  \begin{bmatrix}
    \dot{x}(t)\\
    \dot{I}(t)
  \end{bmatrix}=\begin{bmatrix}
   SW - \beta e_1 e_\ell ^T & k e_1\\
-\theta e^T_\ell   &0
  \end{bmatrix} \begin{bmatrix}
    x(t)\\
    I(t)
  \end{bmatrix}+\begin{bmatrix}
    \frac{\beta}{\theta}\\
    1
  \end{bmatrix}\mu
\end{equation}
where we can immediately recognize the $R$ matrix involved in the Lyapunov equation \eqref{eq:Lyapunov}.

\subsection*{Example - Gene expression network}

\begin{figure}[H]
  \centering
  \includegraphics[width=0.8\textwidth]{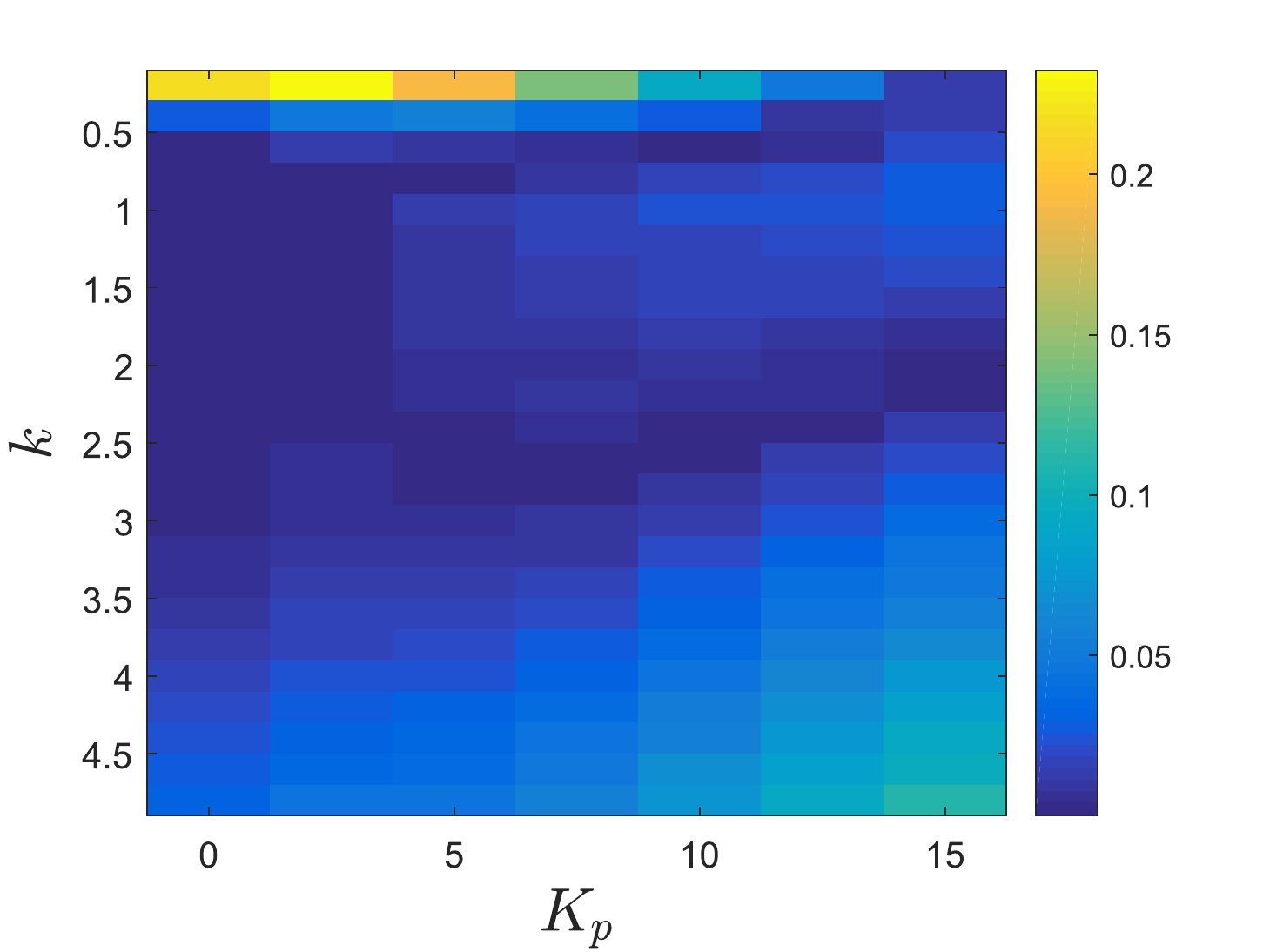}
  \caption{Absolute value of the relative error between the exact stationary variance of the protein copy number and the approximate formula \eqref{eq:vpi_gene} when the gene expression network is controlled with the antithetic integral controller \eqref{eq:AIC} and an ON/OFF proportional controller.}\label{fig:Gene_Prop_RE}
\end{figure}

We present here the results obtained for the gene expression  network \eqref{eq:gene_expression} using the two negative feedback actions. In particular, we will numerically verify the validity of the formula \eqref{eq:variance_ge} and study the influence of the controller parameters on various properties of the closed-loop network. The matrix $R$ is given in this case by
\begin{equation}
  R=\begin{bmatrix}
    -\gamma_r & -\beta & k\\
    k_p & -\gamma_p & 0\\
    0 & -\theta & 0
  \end{bmatrix}.
\end{equation}
It can be shown that the above matrix is Hurwitz stable (i.e. all its eigenvalues are located in the open left half-plane) if and only if the parameters $k,\beta>0$ satisfy the inequality
\begin{equation}\label{eq:RH_gene}
  1 - \dfrac{k \theta k_p}{  \gamma_r \gamma_p (\gamma_r +\gamma_p) }+  \dfrac{ \beta k_p }{ \gamma_r \gamma_p}>0.
\end{equation}
Hence, given $k>0$, the matrix $R$ will be Hurwitz stable for any sufficiently large $\beta>0$ illustrating the stabilizing effect of the proportional action. When the above condition is met, then the closed-loop stationary variance $\V_\pi^{PI}(X_2)$ of the protein copy number is approximately given by the expression
\begin{equation}\label{eq:vpi_gene}
\V_\pi^{PI}(X_2)\approx
\Sigma_{22}= \dfrac{\mu}{\theta}\left[ \dfrac{1 +  \dfrac{k_p}{ \gamma_r +\gamma_p } + \dfrac{k k_p}{ \gamma_r \gamma_p }  +\dfrac{ \beta k_p }{ \gamma_r(  \gamma_r +\gamma_p) }  }{1 - \dfrac{k \theta k_p}{  \gamma_r \gamma_p (\gamma_r +\gamma_p) }+  \dfrac{ \beta k_p }{ \gamma_r \gamma_p}    } \right].
\end{equation}
For any fixed $k>0$ such that \eqref{eq:RH_gene} is satisfied, the closed-loop steady-state variance is a monotonically decreasing function of $\beta$. As a consequence, there will exist a $\beta_c>0$ such that
\begin{equation}
  \Sigma_{22}<\dfrac{\mu}{\theta}\left(1+\dfrac{k_p}{\gamma_r+\gamma_p}\right)
\end{equation}
for all $\beta>\beta_c$.
In particular, when $\beta\to\infty$, then we have that
\begin{equation}
   \Sigma_{22}\to\dfrac{\mu}{\theta}\dfrac{\gamma_p}{\gamma_r+\gamma_p}<\dfrac{\mu}{\theta}.
\end{equation}
We now analyze the results obtained with the antithetic integral controller combined with an OF/OFF proportional feedback. The first step is the numerical comparison of the approximate formula \eqref{eq:vpi_gene} with the stationary variance computed using $10^6$ SSA simulations with the parameters $k_p=2$, $\gamma_r=2$, $\gamma_p=7$, $\mu=10$, $\theta=2$ and $\eta=100$. The absolute value of the relative error between the exact and the approximate stationary variance of the protein copy number for several values for the gains $k$ and $K_p$ is depicted in Figure \ref{fig:Gene_Prop_RE}. We can observe there that the relative error is less than 15\% except when $k$ is very small where the relative error is much larger. However, in this latter case, the mean trajectories do not have time to converge to their steady state value and, therefore, what is depicted in the figure for this value is not very meaningful. In spite of that, we can observe that the approximation is reasonably accurate.

\begin{figure}[H]
  \centering
  \includegraphics[width=0.8\textwidth]{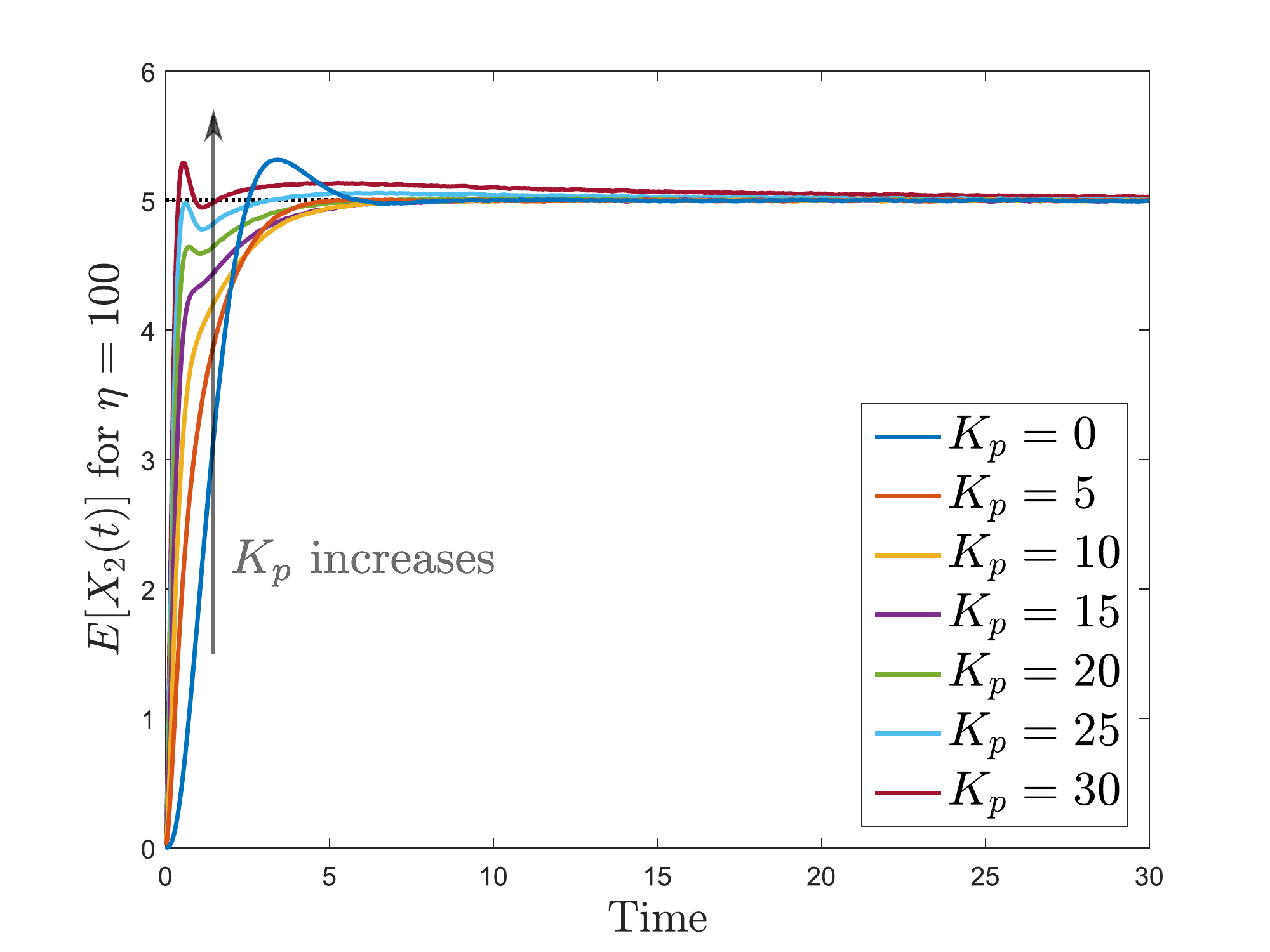}
  \caption{Mean trajectories for the protein copy number when the gene expression network is controlled with the antithetic integral controller \eqref{eq:AIC} with $k=3$ and an ON/OFF proportional controller. The set-point value is indicated as a black dotted line.}\label{fig:Gene_Prop_E}
\end{figure}

\begin{figure}[H]
  \centering
  \includegraphics[width=0.8\textwidth]{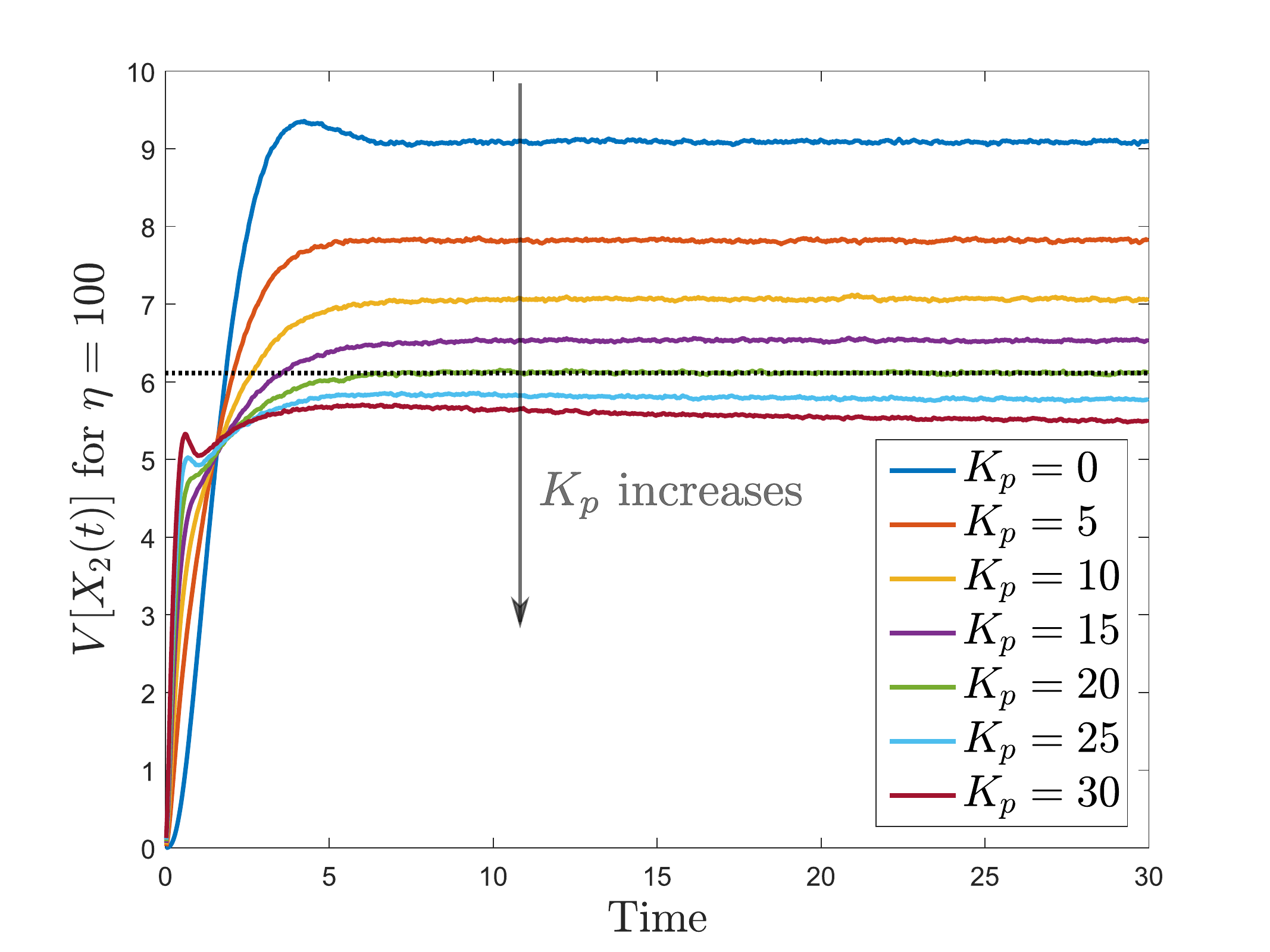}
  \caption{Variance trajectories for the protein copy number when the gene expression network is controlled with the antithetic integral controller \eqref{eq:AIC} with $k=3$ and an ON/OFF proportional controller. The stationary constitutive variance is depicted in black dotted line.}\label{fig:Gene_Prop_V}
\end{figure}

\begin{figure}[H]
  \centering
  \includegraphics[width=0.8\textwidth]{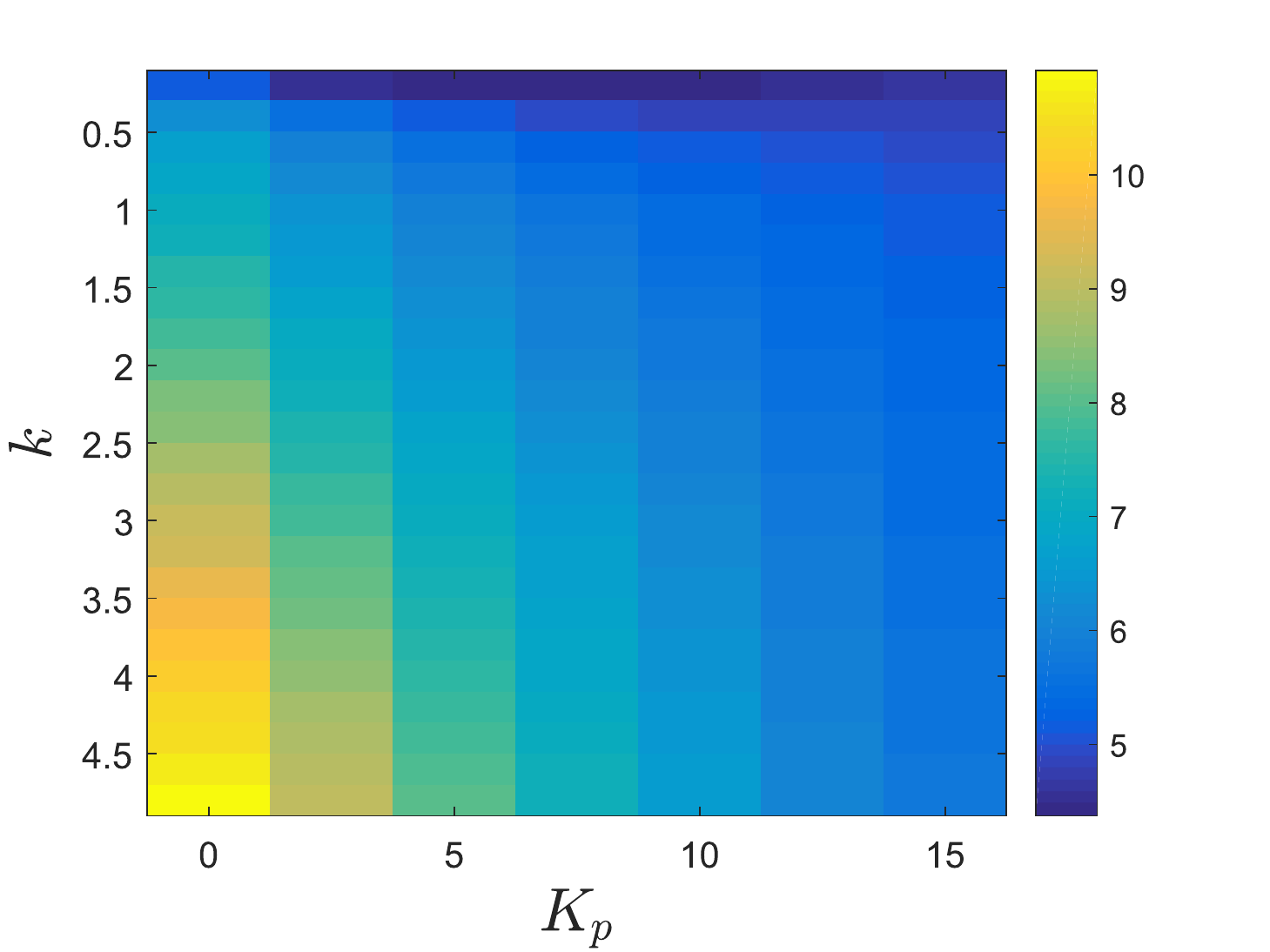}
  \caption{Stationary variance for the protein copy number when the gene expression network is controlled with the antithetic integral controller \eqref{eq:AIC} and an ON/OFF proportional controller.}\label{fig:Gene_Prop_VS}
\end{figure}

\begin{figure}[H]
  \centering
  \includegraphics[width=0.8\textwidth]{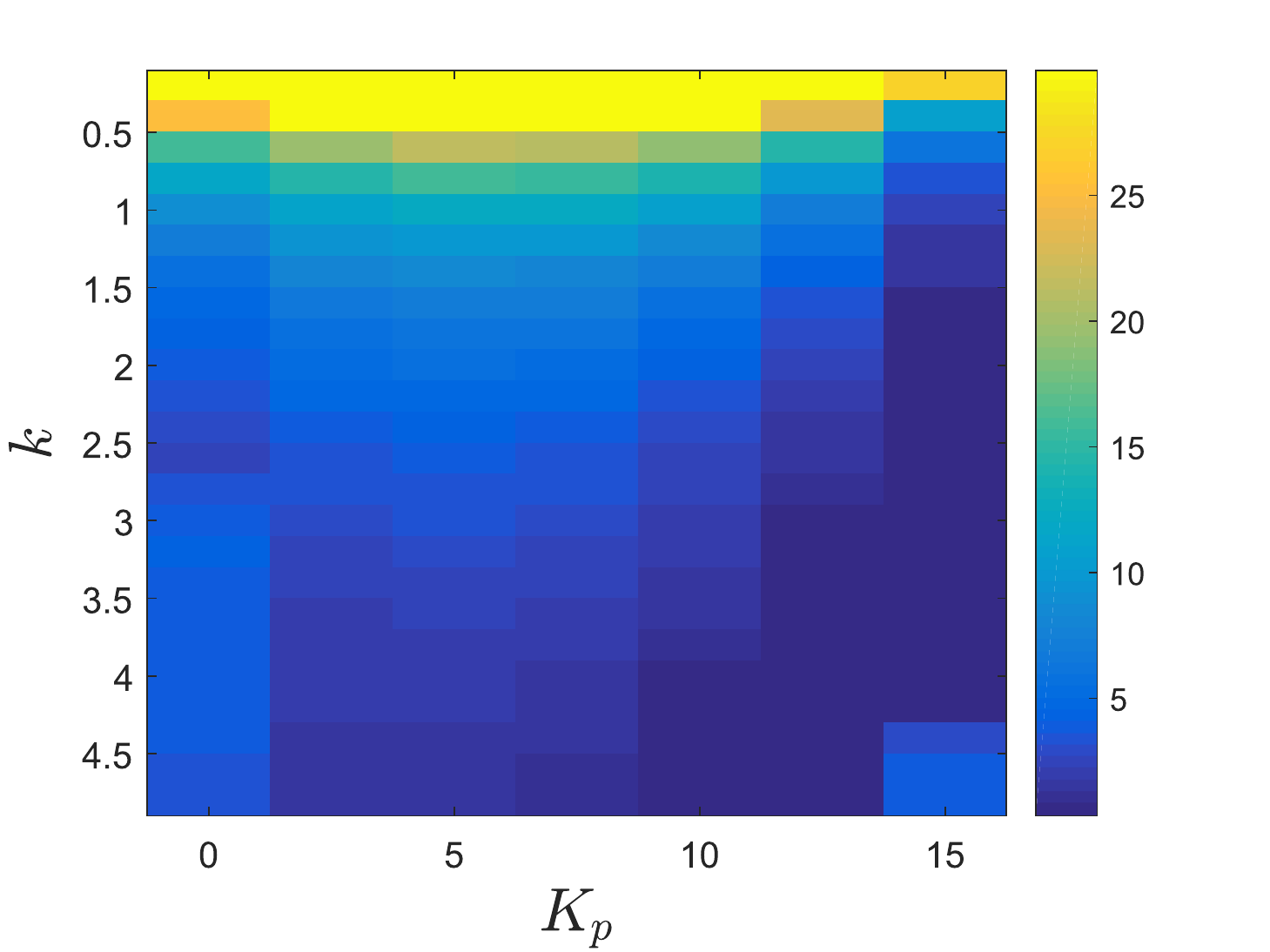}
  \caption{Settling-time for the mean trajectories for the protein copy number when the gene expression network is controlled with the antithetic integral controller \eqref{eq:AIC} and an ON/OFF proportional controller.}\label{fig:Gene_Prop_ST}
\end{figure}

We now look at the performance of the antithetic integral controller combined with an OF/OFF proportional feedback. Figure \ref{fig:Gene_Prop_E} depicts the trajectories of the mean protein copy number while Figure \ref{fig:Gene_Prop_V} depicts the trajectories of the variance of the protein copy number, both in the case where $k=3$. Regarding the mean copy number, we can observe that while at the beginning increasing $K_p$ seems to improve the transient phase, then the dynamics gets more and more abrupt at the start of the transient phase as the gain $K_p$ continues to increase and gets slower and slower at the end of the transient phase, making the means very slow to converge to their set-point. On the other hand, we can see that the stationary variance seems to be a decreasing function of the gain $K_p$. More interestingly, when the gain $K_p$ exceeds 20, the stationary variance becomes smaller than its constitutive value. Figure \ref{fig:Gene_Prop_VS} helps at establishing the influence of the gains $k$ and $K_p$ onto the stationary variance of the protein copy number. We can see that, for any $k$, increasing $K_p$ reduces the stationary variance while for any $K_p$, reducing $k$ reduces the variance, as predicted by the approximate formula \eqref{eq:vpi_gene}. Hence, a suitable choice would be to pick $k$ small and $K_p$ large. We now compare this choice for the parameters with the one that would lead to a small settling-time for the mean dynamics; see Figure \ref{fig:Gene_Prop_ST}. We immediately see that a small $k$ is not an option if one wants to have fast mean dynamics. A sweet spot in this case would be around the right-bottom corner where the settling-time is the smallest. Interestingly, the variance is still at a quite low level even if sometimes higher than the constitutive value.

\begin{figure}[H]
  \centering
  \includegraphics[width=0.8\textwidth]{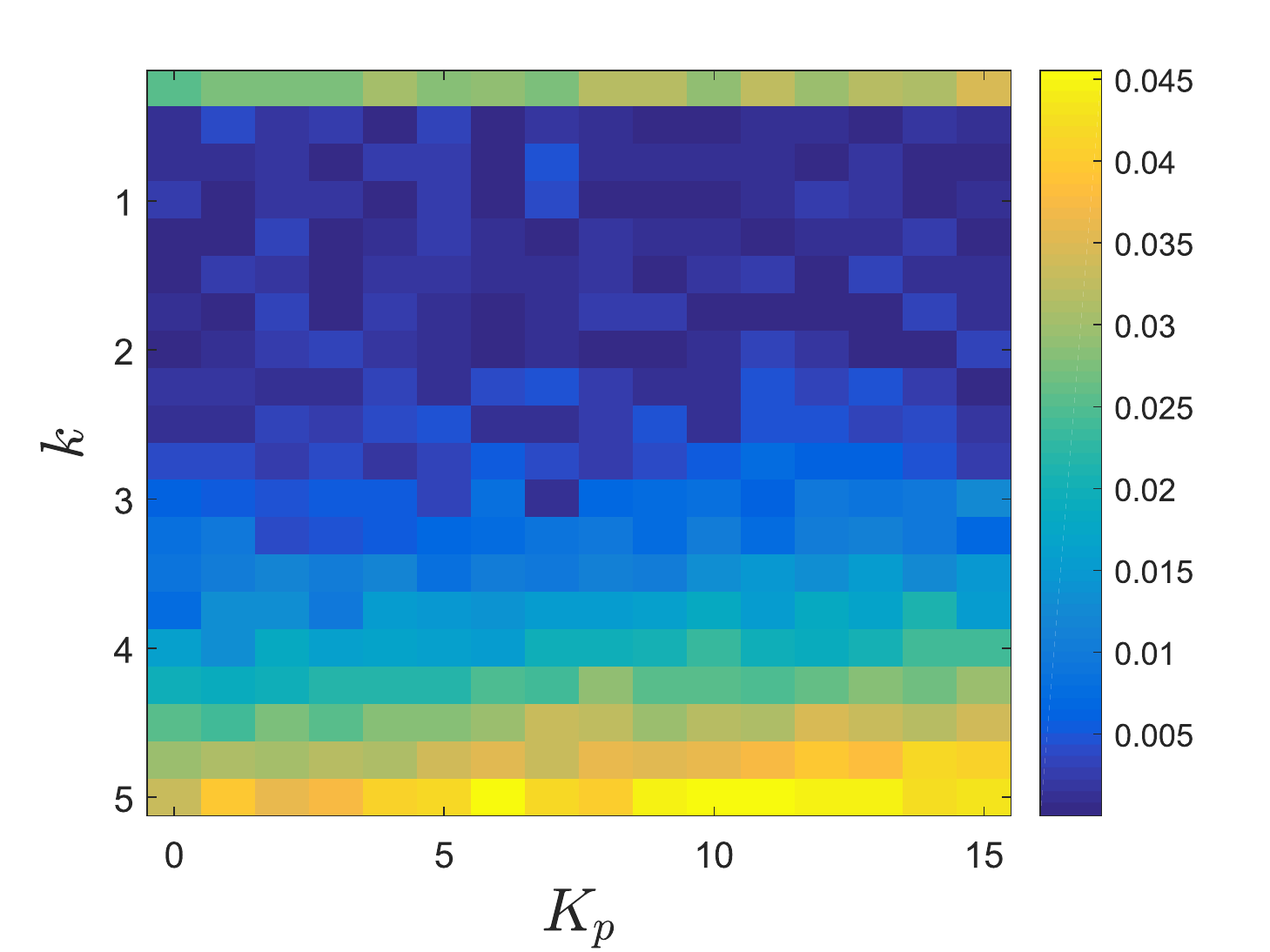}
    \caption{Absolute value of the relative error between the exact stationary variance of the protein copy number and the approximate formula \eqref{eq:vpi_gene} when the gene expression network is controlled with the antithetic integral controller \eqref{eq:AIC} and a Hill controller.}\label{fig:Gene_Hill_RE}
\end{figure}

We now perform the same analysis for the antithetic integral controller combined with the Hill feedback and first verify the accuracy of the approximate formula \eqref{eq:vpi_gene}. We can observe in Figure \ref{fig:Gene_Hill_RE} that the formula is very accurate in this case. To explain this, it is important to note that the gains $K_p$ in both controllers are not directly comparable, only the values for the parameter $\beta$ are. For identical $K_p$'s, the value of $\beta$ for the ON/OFF proportional feedback is much larger than for the Hill feedback (see Figure \ref{fig:Gene_Prop_Beta} and Figure \ref{fig:Gene_Hill_Beta} in the appendix). The Figure \ref{fig:Gene_Prop_RE} and Figure \ref{fig:Gene_Hill_RE} all together simply say that the formula is very accurate when $\beta$ is small.

\begin{figure}[H]
  \centering
  \includegraphics[width=0.8\textwidth]{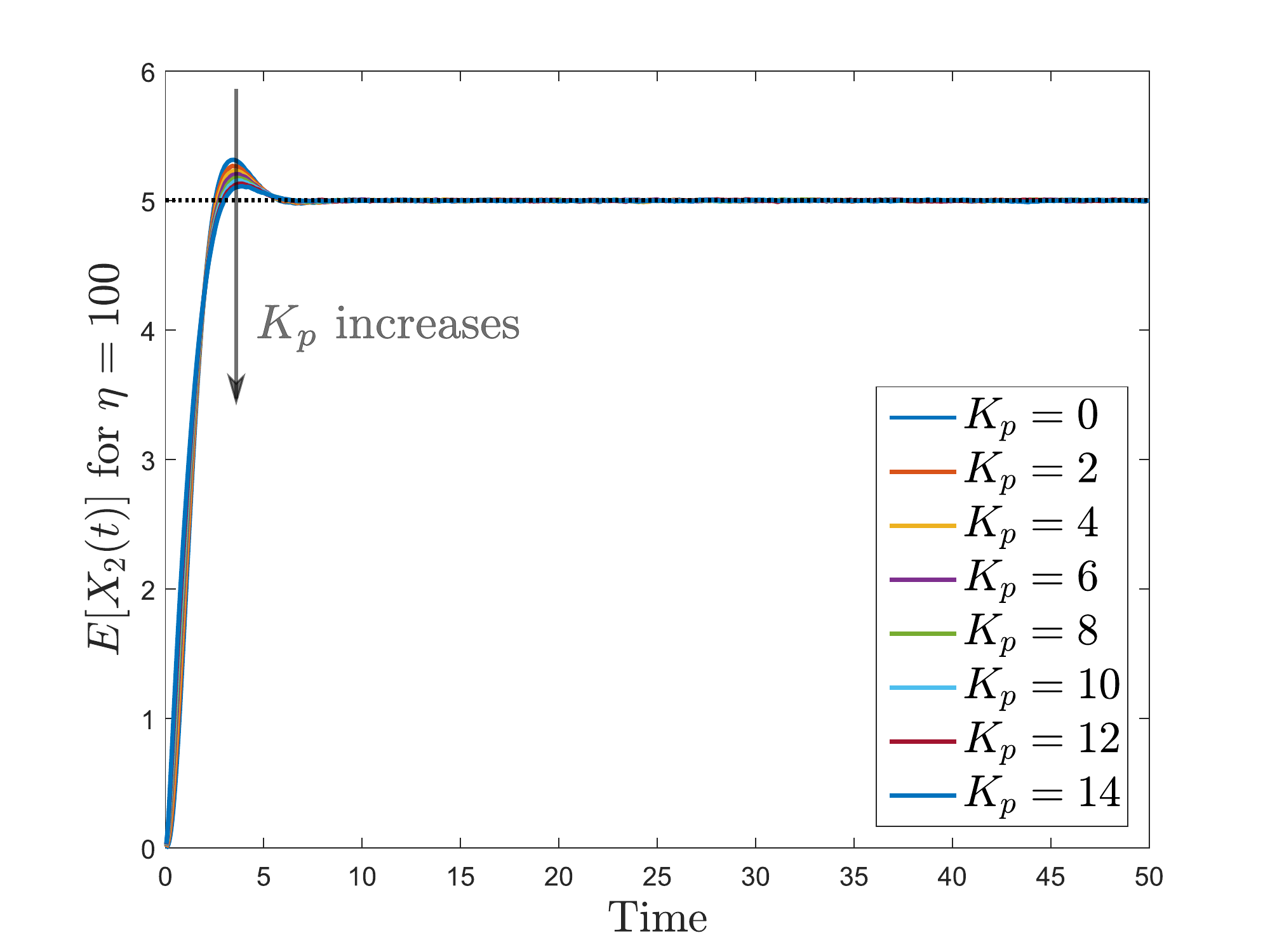}
   \caption{Mean trajectories for the protein copy number when the gene expression network is controlled with the antithetic integral controller \eqref{eq:AIC} with $k=3$ and a Hill controller. The set-point value is indicated as a black dotted line.}\label{fig:Gene_Hill_E}
\end{figure}

\begin{figure}[H]
  \centering
  \includegraphics[width=0.8\textwidth]{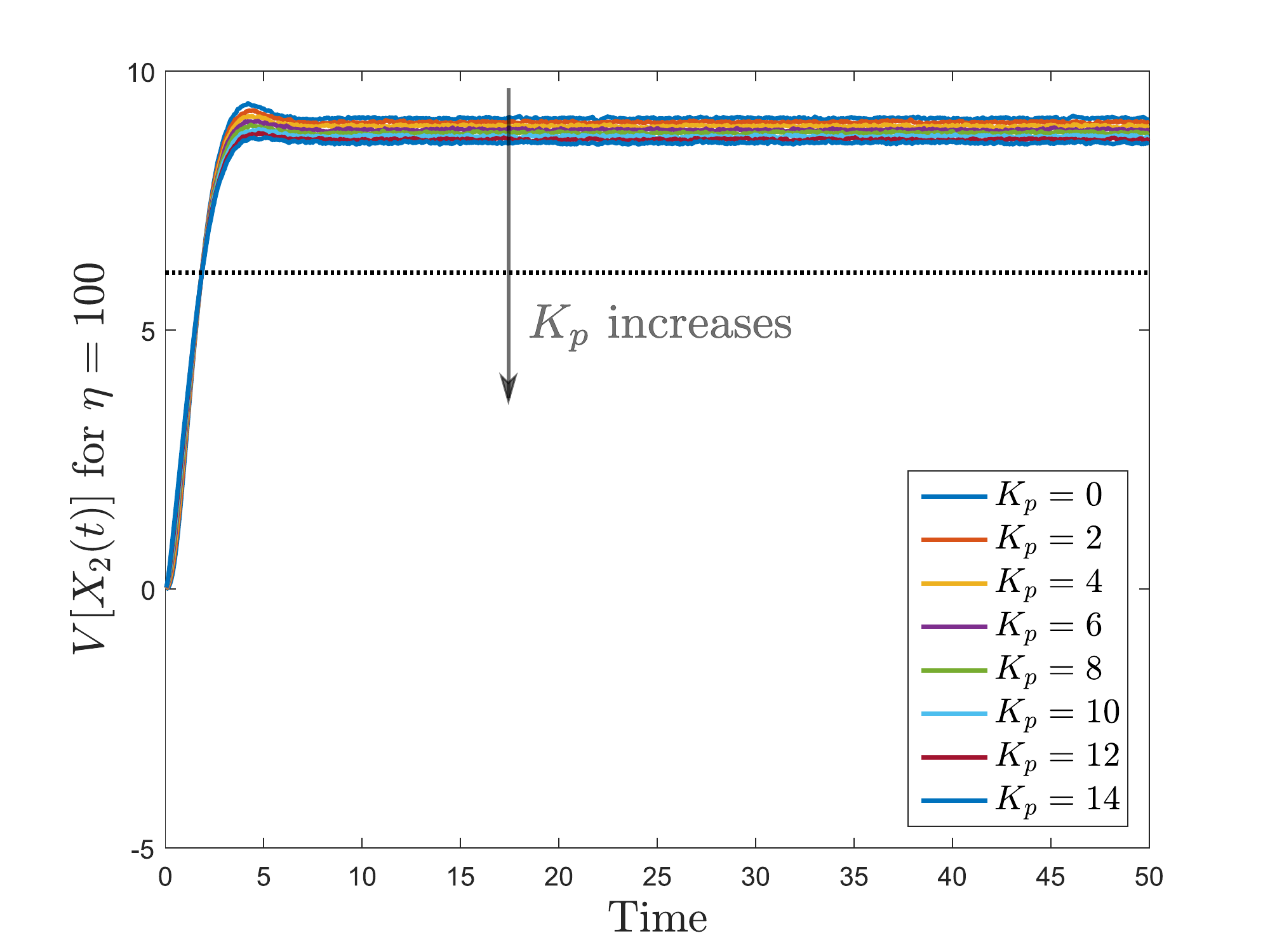}
  \caption{Variance trajectories for the protein copy number when the gene expression network is controlled with the antithetic integral controller \eqref{eq:AIC} with $k=3$ and a Hill controller. The stationary constitutive variance is depicted in black dotted line.}\label{fig:Gene_Hill_V}
\end{figure}

\begin{figure}[H]
  \centering
  \includegraphics[width=0.8\textwidth]{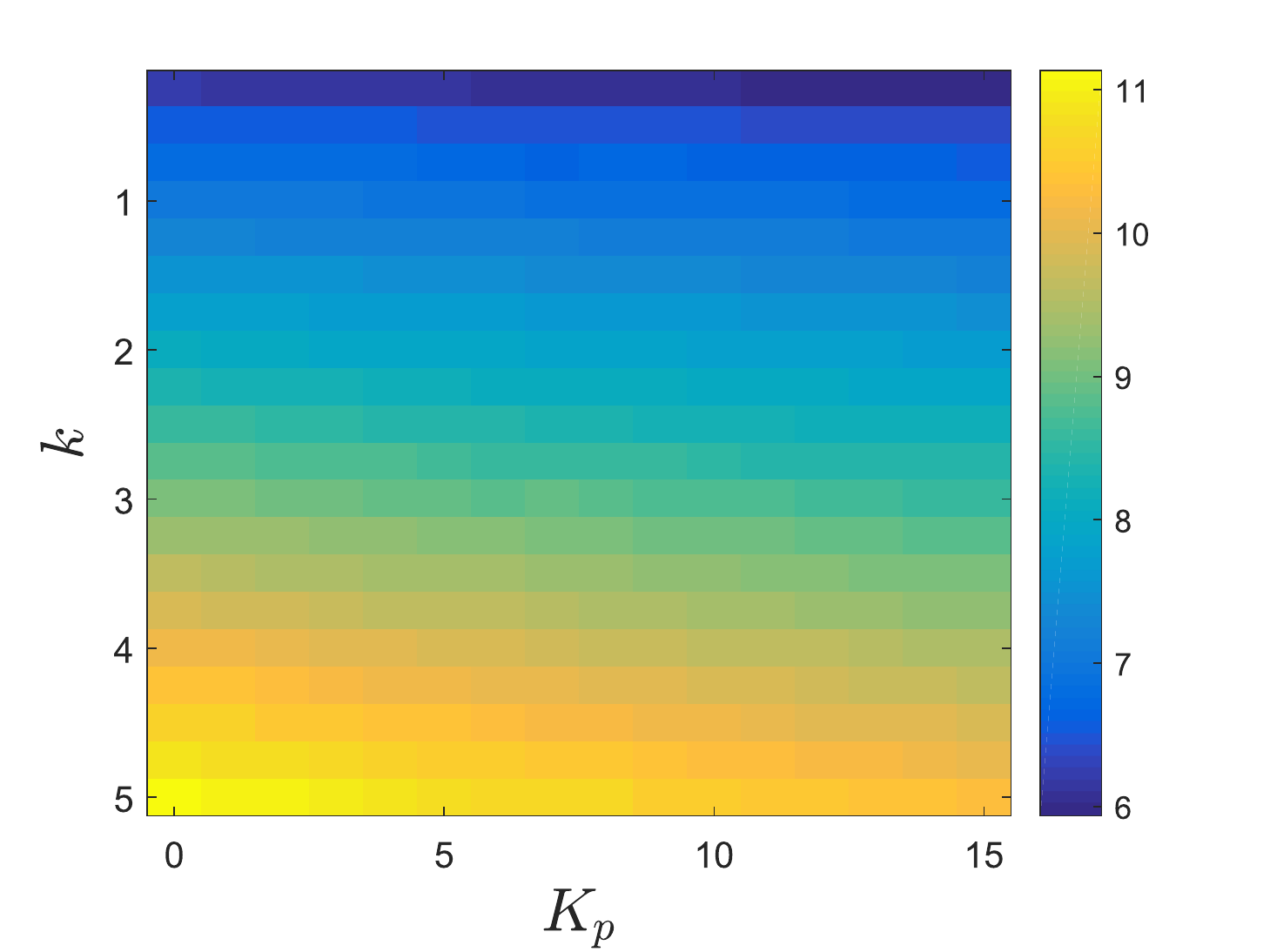}
 \caption{Stationary variance for the protein copy number when the gene expression network is controlled with the antithetic integral controller \eqref{eq:AIC} and a Hill controller.}\label{fig:Gene_Hill_VS}
\end{figure}

\begin{figure}[H]
  \centering
  \includegraphics[width=0.8\textwidth]{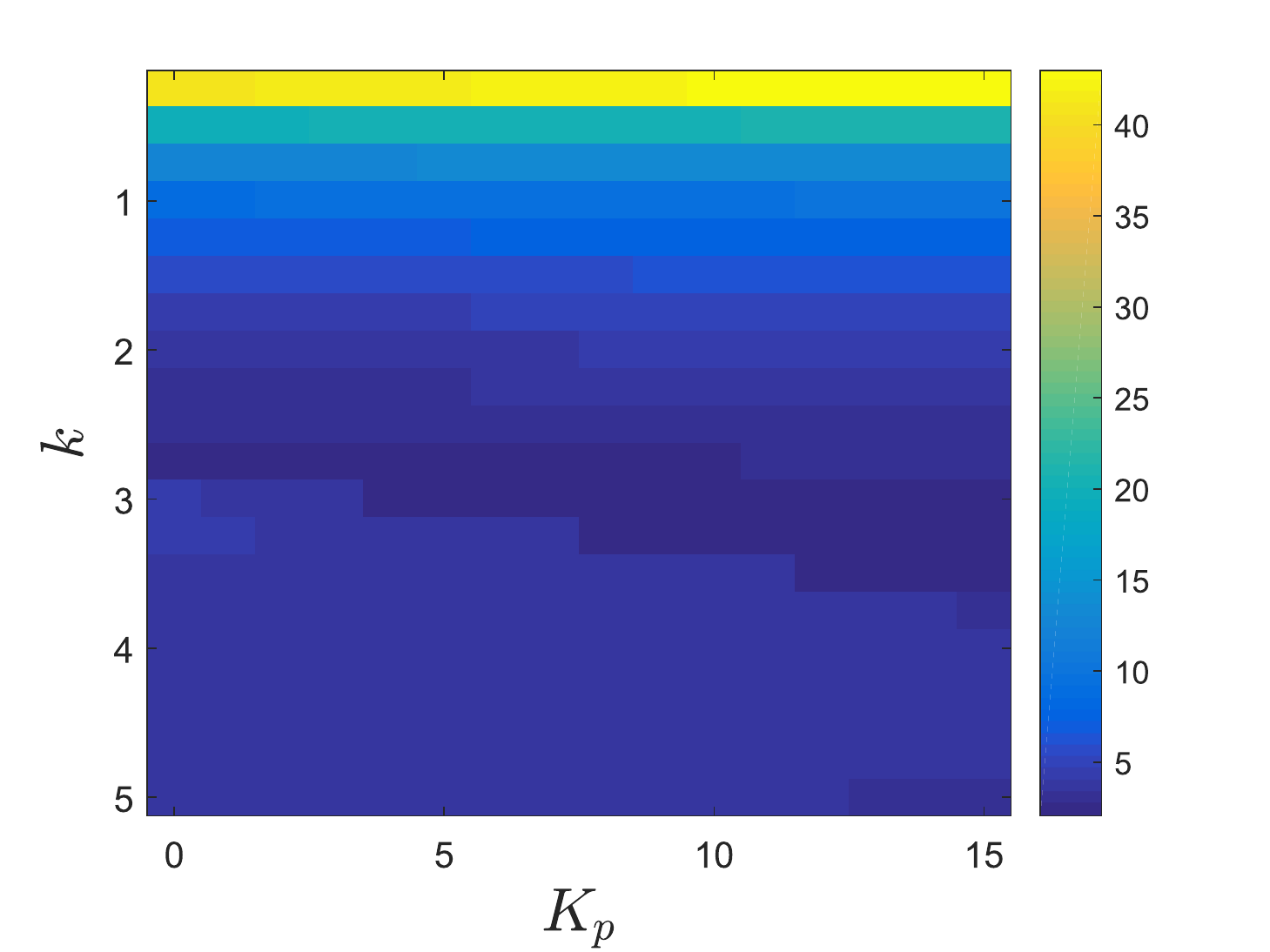}
\caption{Settling-time for the mean trajectories for the protein copy number when the gene expression network is controlled with the antithetic integral controller \eqref{eq:AIC} and a Hill controller.}\label{fig:Gene_Hill_ST}
\end{figure}

We now look at the performance of the antithetic integral controller combined with a Hill feedback.  Similarly to as previously, Figure \ref{fig:Gene_Hill_E} depicts the trajectories of the mean protein copy number while Figure \ref{fig:Gene_Hill_V} depicts the trajectories of the variance of the protein copy number, both in the case where $k=3$. Regarding the mean copy number, we can observe than the dynamics are much more homogeneous than in the previous case and that increasing $K_p$ reduces the overshoot and, hence, the settling-time. This can again be explained by the fact that $\beta$ is much smaller in this case. Similarly, the spread of the variances is much tighter than when using the other negative feedback again because of the fact that $\beta$ is small in this case. This homogeneity is well illustrated in Figure \ref{fig:Gene_Hill_VS} and Figure \ref{fig:Gene_Hill_ST} where we conclude on the existence of a clear tradeoff between settling-time and stationary variance.

As can been seen in Figure \ref{fig:Gene_Prop_E} and Figure \ref{fig:Gene_Hill_E}, the mean dynamics are quite different and it would be interesting to explain this difference in terms of control theoretic ideas. A first explanation lies in the sensitivity of the parameter $\beta$ in terms of the feedback strength $K_p$. In the case of the ON/OFF proportional feedback, this sensitivity is quite high whereas it is very low in the case of the Hill feedback (see Figure \ref{fig:Gene_Prop_Beta} and Figure \ref{fig:Gene_Hill_Beta} in the appendix). This gives an explanation on why the mean trajectories are very different in the case of the ON/OFF proportional feedback for different values of $K_p$ while the mean trajectories are very close to each other in the case of the Hill feedback. A second explanation lies in the type of feedback in use. Indeed, the ON/OFF proportional feedback is an error-feedback and, when combined with the antithetic integral controller, may introduce a stable zero in the mean dynamics. On the other hand, the Hill feedback is an output-feedback that does not seem to introduce such a zero. When increasing the negative feedback gain $K_p$, this zero moves towards the origin. Once very close to the origin, this zero will have an action in the closed-loop mean dynamics that is very close to a derivative action, leading then to abrupt initial transient dynamics. A theoretical basis for this discussion is developed in more details in the SI.

\subsection*{Example - Gene expression network with protein maturation}

\begin{figure}[H]
  \centering
  \includegraphics[width=0.8\textwidth]{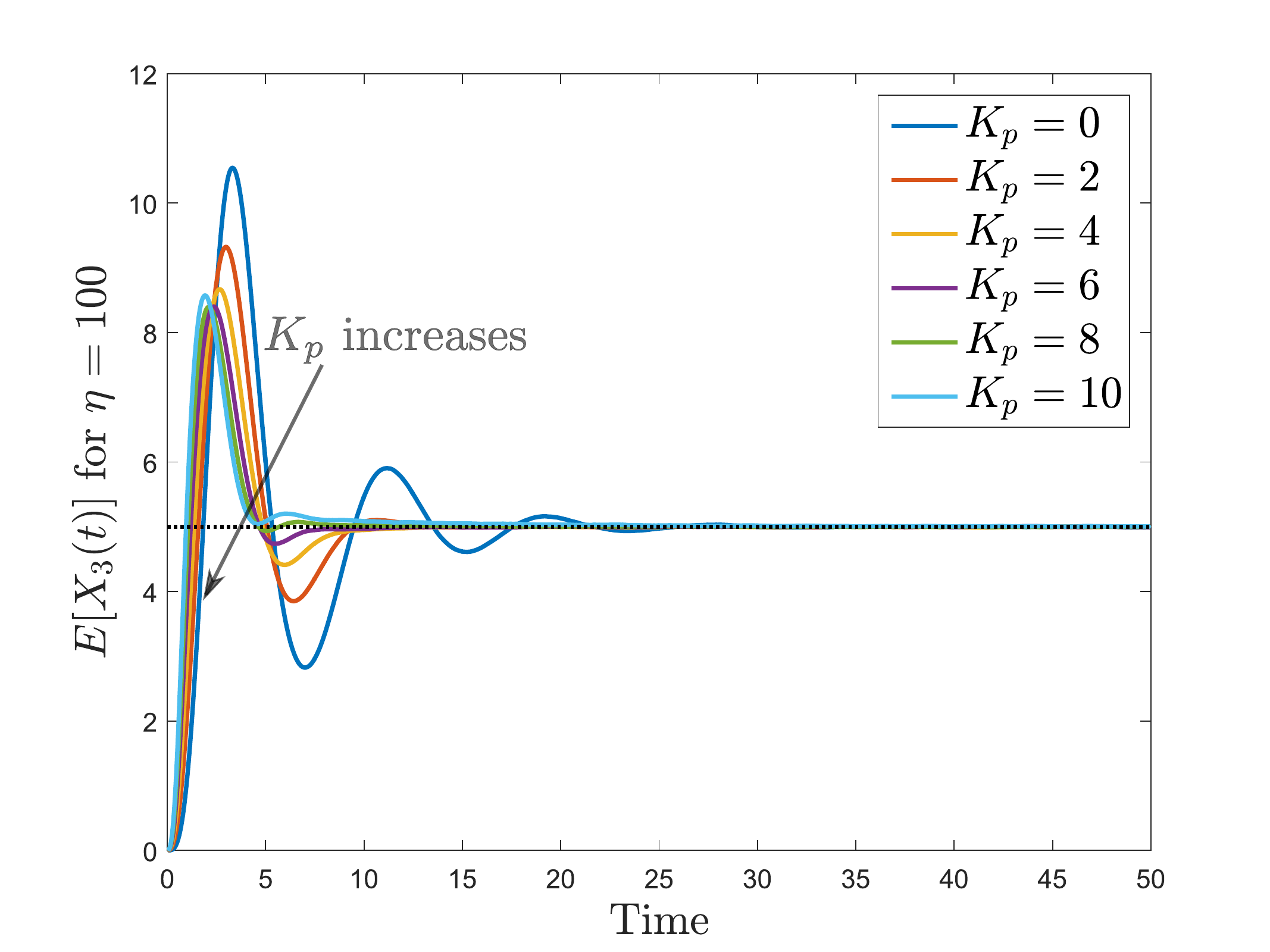}
    \caption{Mean trajectories for the mature protein copy number when the gene expression network with protein maturation is controlled with the antithetic integral controller \eqref{eq:AIC} with $k=3$ and an ON/OFF proportional controller. The set-point value is indicated as a black dotted line.}\label{fig:Mat_Prop_E}
\end{figure}

\begin{figure}[H]
  \centering
  \includegraphics[width=0.8\textwidth]{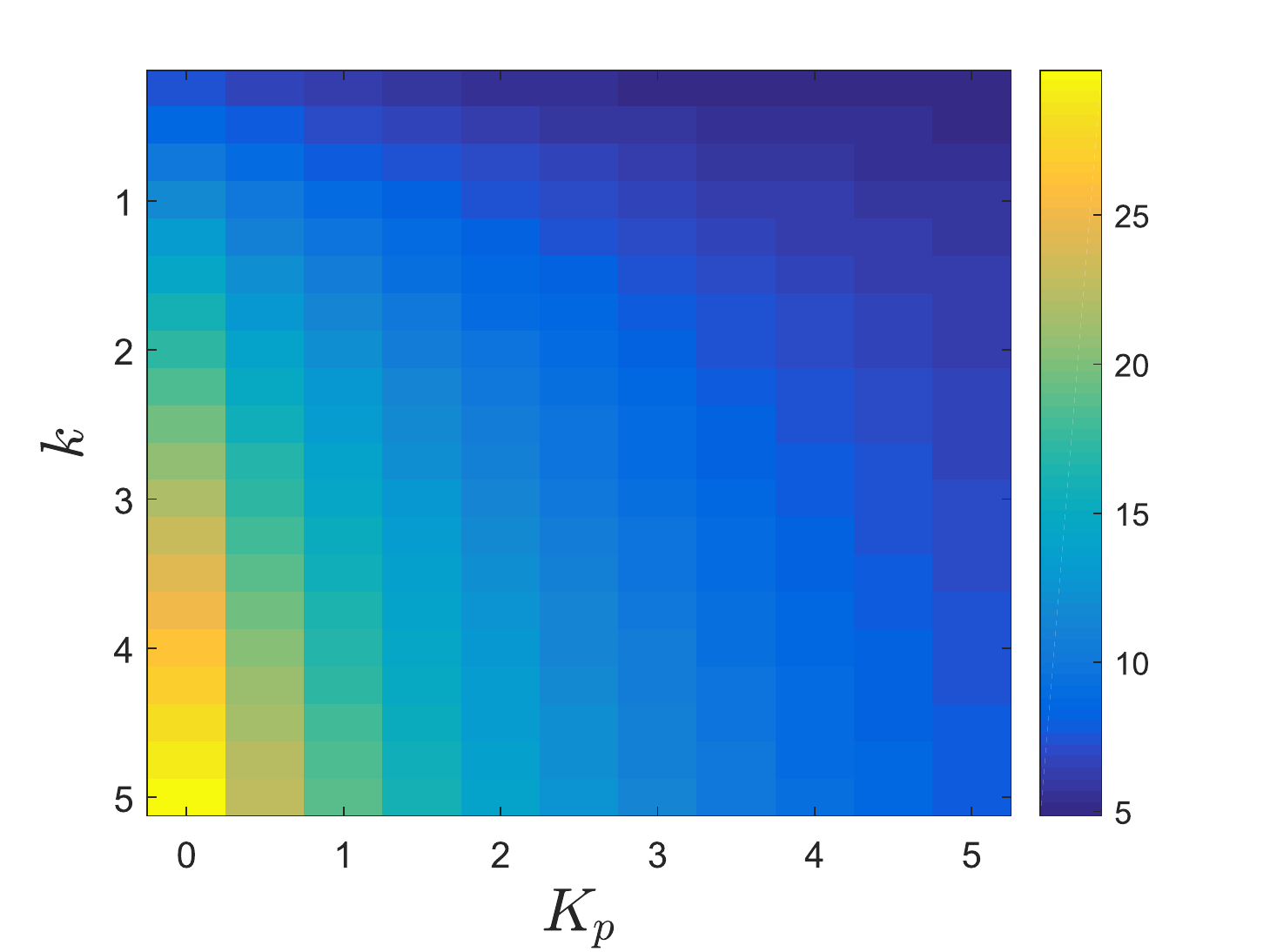}
    \caption{Stationary variance for the mature protein copy number when the gene expression network with protein maturation is controlled with the antithetic integral controller \eqref{eq:AIC} and an ON/OFF proportional controller.}\label{fig:Mat_Prop_VS}
\end{figure}

\begin{figure}[H]
  \centering
  \includegraphics[width=0.8\textwidth]{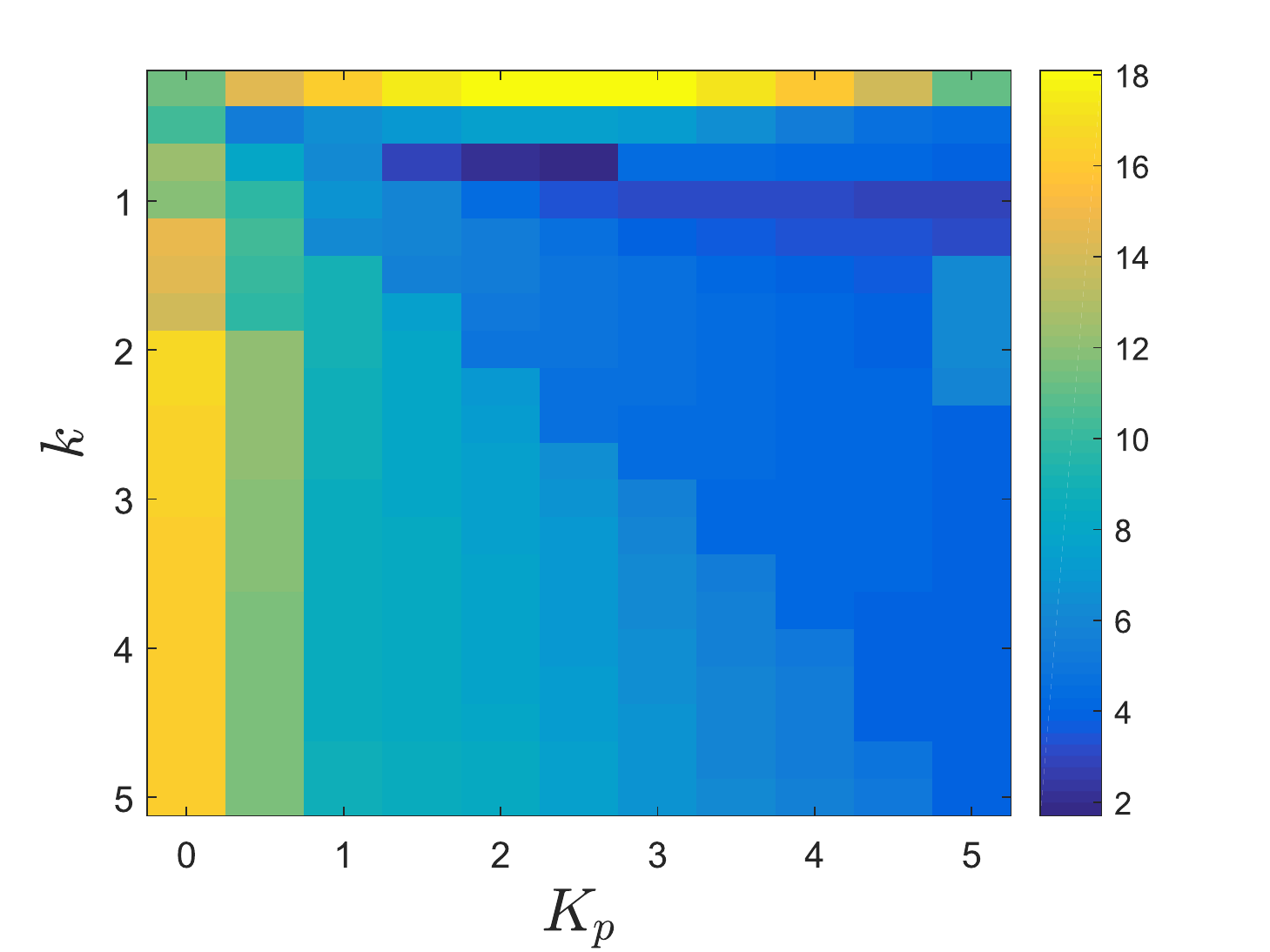}
    \caption{Settling-time for the mean trajectories for the mature protein copy number when the gene expression network with protein maturation is controlled with the antithetic integral controller \eqref{eq:AIC} and an ON/OFF proportional controller.}\label{fig:Mat_Prop_ST}
\end{figure}

The results obtained in the previous section clearly only hold for the gene expression network and it would be quite hasty to directly generalize those results to more complex unimolecular networks. This hence motivates the consideration of a slightly more complicated example, namely, the gene expression network involving a protein maturation reaction given by
  \begin{equation}
  \begin{array}{c}
  \phib\rarrow{k_r}\X{1},\ \X{1}\rarrow{k_p}\X{1}+\X{2}, \X{1}\rarrow{\gamma_r}\phib, \X{2}\rarrow{\gamma_p}\phib\\
      \X{2}\rarrow{k_p'}\X{3}, \X{3}\rarrow{\gamma_p'}\phib
  \end{array}
\end{equation}
where, as before, $\X{1}$ denotes mRNA, $\X{2}$ denotes protein and, now, $\X{3}$ denotes the mature protein. In this case, the goal is to control the average mature protein copy number by, again, acting at a transcriptional level. As this network is still unimolecular, the proposed framework remains valid. In particular, the matrix $R$ is given by
\begin{equation}
  R=\begin{bmatrix}
    -\gamma_r & 0 & -\beta & k\\
    k_p & -(\gamma_p+k_p') & 0 & 0\\
    0 & k_p' & -\gamma_p' & 0\\
    0 & 0 & -\theta & 0
  \end{bmatrix}
\end{equation}
and is Hurwitz stable provided that the two following conditions are satisfied
\begin{equation}
  \beta<\dfrac{1}{k_pk_p'}\left((\gamma_r + \gamma_p + \gamma_p' + k_p')(\gamma_r\gamma_p +\gamma_r\gamma_p. + \gamma_p\gamma_p' + \gamma_rk_p' + \gamma_p'k_p')-\gamma_r\gamma_p'(\gamma_p+k_p')\right)
\end{equation}
and
\begin{equation}
  k_p^2k_p'^2\beta^2+\sigma_1\beta+\sigma_0<0
\end{equation}
where
\begin{equation}
  \begin{array}{rcl}
    \sigma_1&=&- k_pk_p'(\gamma_r + \gamma_p'+ \gamma_p + k_p')(\gamma_r\gamma_p + \gamma_r\gamma_p' + \gamma_p\gamma_p' + \gamma_rk_p' + \gamma_p'k_p')\\
    && +2\gamma_r\gamma_p'k_pk_p'(\gamma_p + k_p'),\\
    \sigma_0&=& - \gamma_r\gamma_p'(\gamma_p + k_p')(\gamma_r + \gamma_p + \gamma_p' + k_p')(\gamma_r\gamma_p + \gamma_r\gamma_p' + \gamma_p\gamma_p' + \gamma_rk_p' + \gamma_p'k_p')\\
    && +\gamma_r^2\gamma_p2^2(\gamma_p + k_p')^2+ kk_pk_p'\theta(\gamma_r + \gamma_p + \gamma_p' + k_p')^2.
  \end{array}
\end{equation}
Considering, for instance, the following parameters $k_p=1$, $\gamma_r=2$, $\gamma_p=1$, $k_p'=3$, $\gamma_p'=1$, $\mu=10$, $\theta=2$ and $\eta=100$, the above conditions reduce to
\begin{equation}
  \beta<30
\end{equation}
and
\begin{equation}
 9\beta^2 - 246\beta + 294k - 720<0.
\end{equation}
The intersection of these conditions yield the stability conditions
\begin{equation}
  k\in(0,49/6) \textnormal{ and } \beta\in\left(\dfrac{41-7\sqrt{49-6k}}{3}, \dfrac{41+7\sqrt{49-6k}}{3}\right)\cap(0,\infty).
\end{equation}
It can be verified that for values on the boundary of at least one of those intervals, the matrix $R$ has eigenvalues on the imaginary axis. Standard calculations on the moments equation show that the open-loop variance is given by
\begin{equation}\label{eq:vpi_mat}
  \V_\pi^{OL}(X_3)=\dfrac{\mu}{\theta}\left(1+k_pk_p'\dfrac{k_p' + \gamma_r + \gamma_p + \gamma_p'}{(\gamma_r + \gamma_p')(\gamma_r + \gamma_p + k_p')(\gamma_p + \gamma_p' + k_p')}\right).
\end{equation}
%
With the numerical values for the parameters previously given, the open-loop variance is approximately equal to $37/6\approx6.1667$. The closed-loop variance, however, is approximately given by
\begin{equation}
  \V_\pi^{PI}(X_3)\approx\Sigma_{33}=\dfrac{\mu}{\theta}\left(\dfrac{\dfrac{\theta}{\mu} \V_\pi^{OL}(X_3)+\dfrac{\zeta_k}{\zeta_d}k+\dfrac{\zeta_\beta}{\zeta_d}\beta+\dfrac{\zeta_{k\beta}}{\zeta_d}k\beta}{1+\dfrac{\xi_k}{\xi_d}k+\dfrac{\xi_\beta}{\xi_d}\beta+\dfrac{\xi_{\beta^2}}{\xi_d}\beta^2}\right)
\end{equation}
where
\begin{equation}
  \begin{array}{rcl}
    \xi_d&=&\gamma_r\gamma_p'(\gamma_r + \gamma_p')(\gamma_p + k_p')(\gamma_r + \gamma_p + k_p')(\gamma_p + \gamma_p' + k_p')\\
    \xi_k&=&-k_pk_p'\theta(\gamma_r + \gamma_p + \gamma_p' + k_p')^2\\
    \xi_\beta&=&k_pk_p'(\gamma_r^2\gamma_p + \gamma_r^2\gamma_p' + \gamma_r^2k_p' + \gamma_r\gamma_p^2 + \gamma_r\gamma_p\gamma_p' + 2\gamma_r\gamma_pk_p' + \gamma_r\gamma_p'^2\\
     &&+ \gamma_r\gamma_p'k_p' + \gamma_rk_p'^2 + \gamma_p^2\gamma_p' + \gamma_p\gamma_p'^2 + 2\gamma_p\gamma_p'k_p' + \gamma_p'^2k_p' + \gamma_p'k_p'^2)\\
    \xi_{\beta^2}&=&-k_p^2k_p'^2
  \end{array}
\end{equation}
and
\begin{equation}
  \begin{array}{rcl}
    \zeta_d&=&\xi_d\\
    \zeta_k&=&k_pk_p'(\gamma_r^2\gamma_p + \gamma_r^2\gamma_p' + \gamma_r^2k_p' + \gamma_r\gamma_p^2 + 2\gamma_r\gamma_p\gamma_p' + 2\gamma_r\gamma_pk_p' + \gamma_r\gamma_p'^2\\
     &&+ 2\gamma_r\gamma_p'k_p' - \theta\gamma_r\gamma_p' + \gamma_rk_p'^2 + \gamma_p^2\gamma_p' + \gamma_p\gamma_p'^2 + 2\gamma_p\gamma_p'k_p' - \theta\gamma_p\gamma_p' \\
     &&+ \gamma_p'^2k_p' - \theta\gamma_p'^2 + \gamma_p'k_p'^2 - \theta\gamma_p'k_p')\\
    \zeta_\beta&=&\gamma_p'k_pk_p'(\gamma_r^2 + \gamma_r\gamma_p + \gamma_rk_p' + \gamma_p'\gamma_r + \gamma_p^2 + 2\gamma_pk_p' + \gamma_p'\gamma_p + k_p'^2 + \gamma_p'k_p')\\
    \zeta_{k\beta}&=&-k_p^2k_p'^2.
  \end{array}
\end{equation}
An expression that is more complex than, yet very similar to, the formula \eqref{eq:vpi_gene} obtained for the simple gene expression network. For the considered set of parameter values, the approximated variance is a nonmonotonic function of the parameter $\beta$ as it can be theoretically observed in Figure \ref{fig:Mat_NM} in the appendix. It turns out that this behavior can also be observed in the numerical simulations depicted in Figure \ref{fig:Mat_Prop_V_NM} in the appendix where we can see that the variance exhibits this monotonic behavior. However, it should also be pointed out that this increase is accompanied with the emergence of an tracking error for the mean dynamics (see Figure \ref{fig:Mat_Prop_V_NM} in the appendix) and a loss of ergodicity for the overall controlled network as emphasized by diverging mean dynamics for the sensing species (see Figure \ref{fig:Mat_Prop_Z2_NM} in the appendix).  This contrasts with the gene expression case where the variance was a monotonically decreasing function of $\beta$. Regarding the mean dynamics, we can see that increasing $K_p$ and, hence, $\beta$, to reasonable levels improves the settling-time as depicted in Figure \ref{fig:Mat_Prop_E} for the special case of $k=3$. However, this is far from being the general case since the settling-time can exhibit a quite complex behavior for this network (see \ref{fig:Mat_Prop_ST}). The stationary variance depicted in Figure  \ref{fig:Mat_Prop_VS} exhibits here a rather standard and predictive behavior where a small $k$ and a large $K_p$ both lead to its reduction. Similar conclusions can be drawn when the network is controlled with a Hill negative feedback controller; see Figure \ref{fig:Mat_Hill_E}, Figure \ref{fig:Mat_Hill_E}, Figure \ref{fig:Mat_Hill_E} and Figure \ref{fig:Mat_Hill_E} in the appendix.

\begin{figure}[H]
  \centering
  \includegraphics[width=0.8\textwidth]{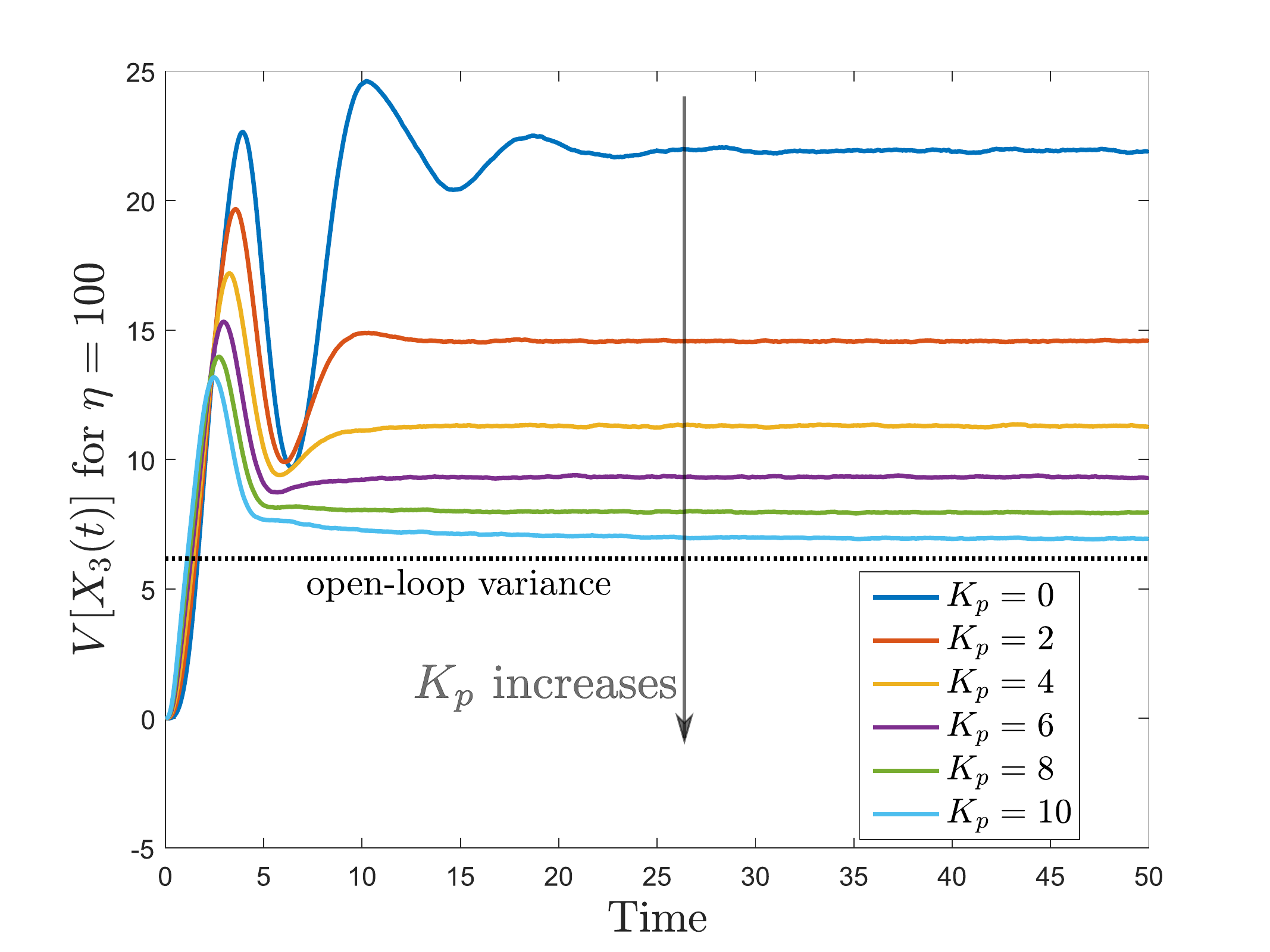}
    \caption{Variance trajectories for the mature protein copy number when the gene expression network with protein maturation is controlled with the antithetic integral controller \eqref{eq:AIC} with $k=3$ and an ON/OFF proportional controller. The stationary constitutive variance is depicted in black dotted line.}\label{fig:Mat_Prop_V}
\end{figure}

\subsection*{Example - Gene expression network with protein dimerization}

\begin{figure}[H]
  \centering
  \includegraphics[width=0.8\textwidth]{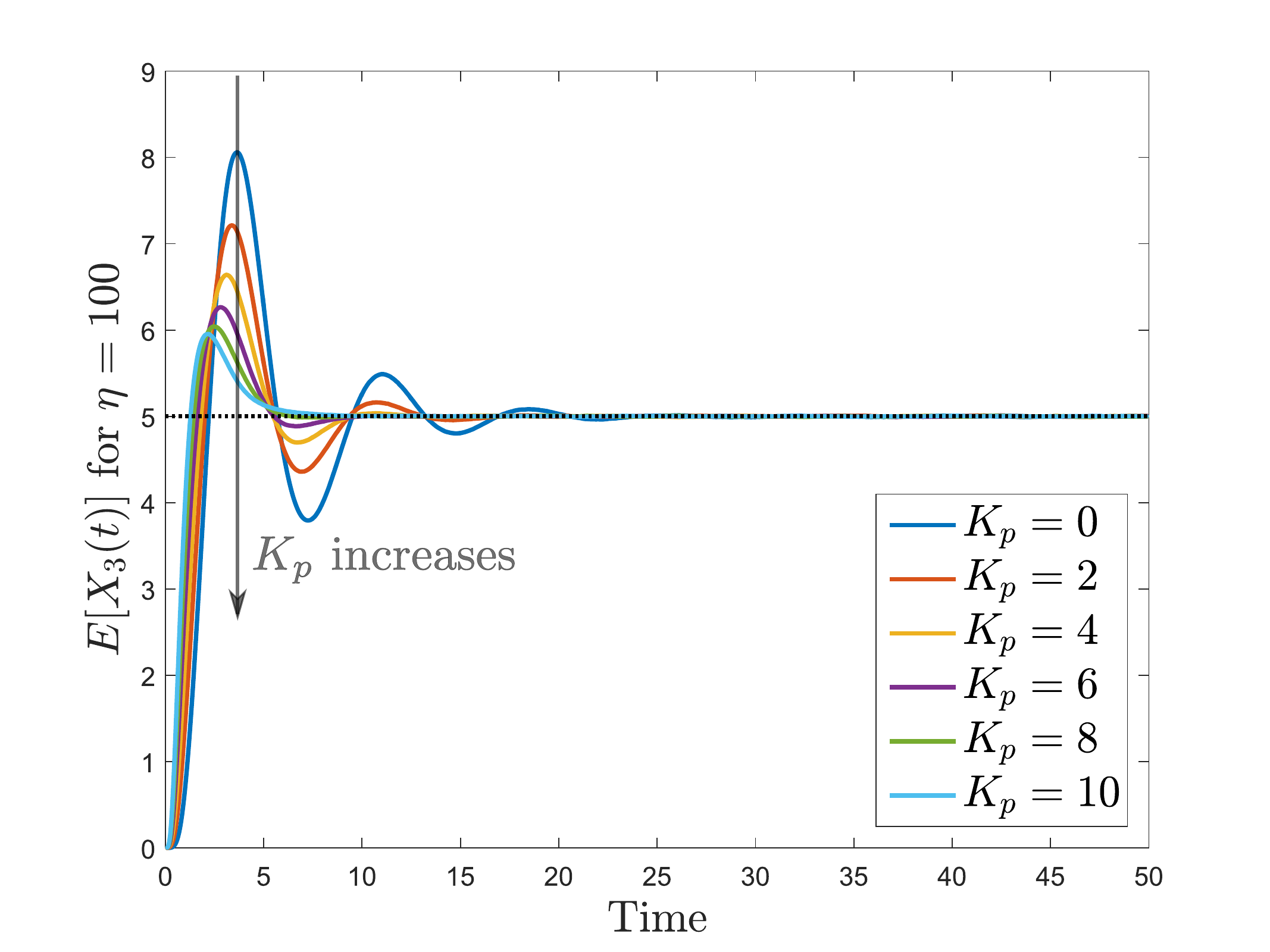}
  \caption{Mean trajectories for the homodimer copy number when the gene expression network with protein dimerization is controlled with the antithetic integral controller \eqref{eq:AIC} with $k=3$ and an ON/OFF proportional controller. The set-point value is indicated as a black dotted line.}\label{fig:Dimer_Prop_E}
\end{figure}

\begin{figure}[H]
  \centering
  \includegraphics[width=0.8\textwidth]{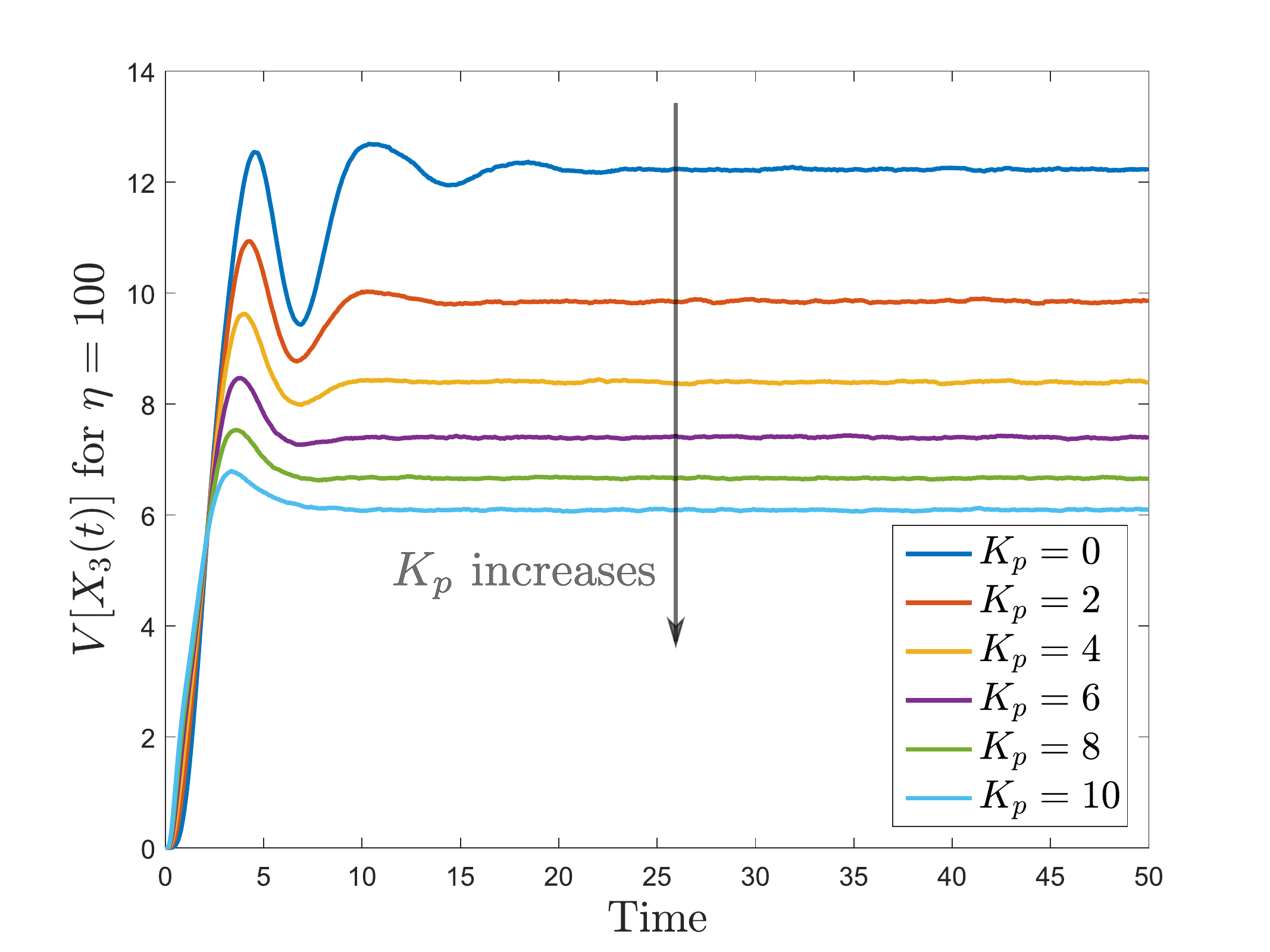}
 \caption{Variance trajectories for the homodimer copy number when the gene expression network with protein dimerization  is controlled with the antithetic integral controller \eqref{eq:AIC} with $k=3$ and an ON/OFF proportional controller. The stationary constitutive variance is depicted in black dotted line.}\label{fig:Dimer_Prop_V}
\end{figure}

\begin{figure}[H]
  \centering
  \includegraphics[width=0.8\textwidth]{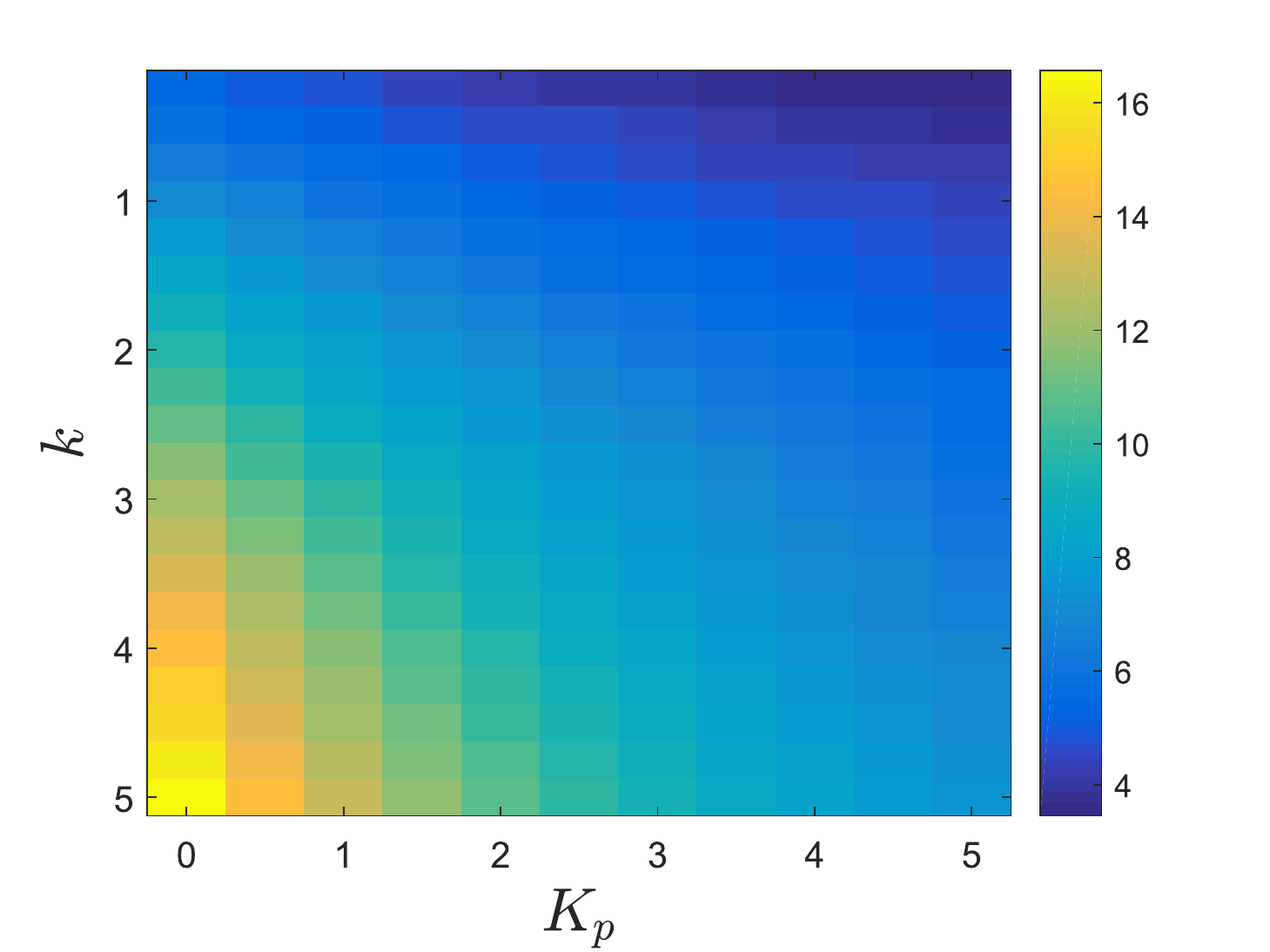}
  \caption{Stationary variance for the homodimer copy number when the gene expression network with protein dimerization  is controlled with the antithetic integral controller \eqref{eq:AIC} and an ON/OFF proportional controller.}\label{fig:Dimer_Prop_VS}
\end{figure}

\begin{figure}[H]
  \centering
  \includegraphics[width=0.8\textwidth]{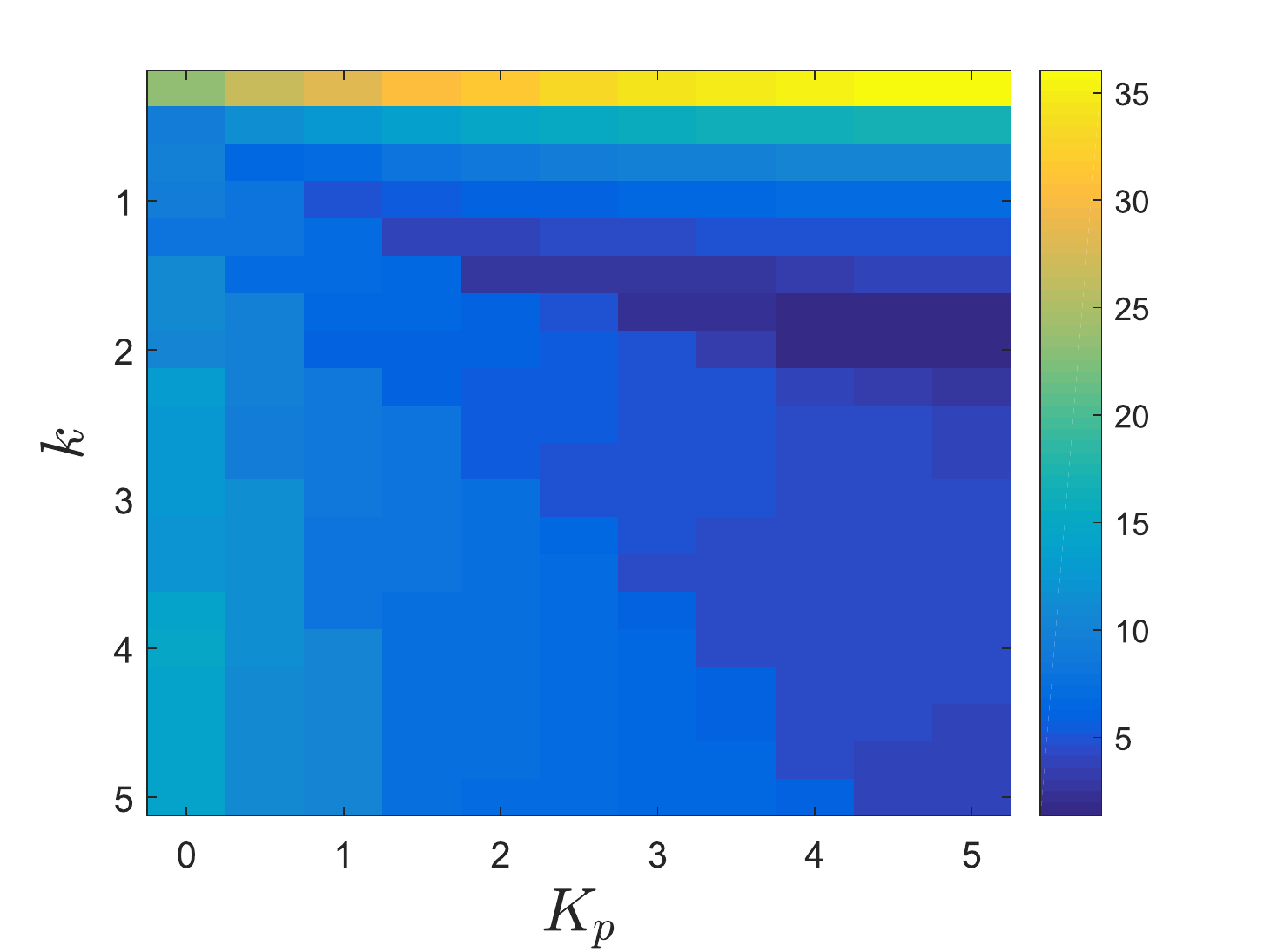}
  \caption{Settling-time for the mean trajectories for the homodimer copy number when the gene expression network with protein dimerization  is controlled with the antithetic integral controller \eqref{eq:AIC} and an ON/OFF proportional controller.}\label{fig:Dimer_Prop_ST}
\end{figure}

The proposed theory is only valid for unimolecular networks but, in spite of that, it is still interesting to see whether similar conclusions could be obtained for a network that is not unimolecular. This motivates the consideration of the following gene expression network with protein dimerization:
  \begin{equation}
  \begin{array}{c}
  \phib\rarrow{k_r}\X{1},\ \X{1}\rarrow{k_p}\X{1}+\X{2}, \X{1}\rarrow{\gamma_r}\phib, \X{2}\rarrow{\gamma_p}\phib\\
  \X{2}+\X{2}\rarrow{k_d}\X{3}, \X{3}\rarrow{\gamma_d} \X{2}+\X{2}, \X{3}\rarrow{\gamma_d'}\phib
  \end{array}
\end{equation}
where, as before, $\X{1}$ denotes mRNA, $\X{2}$ denotes protein but, now, $\X{3}$ denotes a protein homodimer. In this case, the Lyapunov equation \eqref{eq:Lyapunov} is not valid anymore because of the presence of the dimerization reaction but we can still perform stochastic simulations. The considered parameter values are given by $k_p=1$, $\gamma_r=2$, $\gamma_p=1$, $k_d=3$, $\gamma_d=\gamma_d'=1$, $\mu=10$, $\theta=2$ and $\eta=100$. We can see in Figure \ref{fig:Dimer_Prop_E}, Figure \ref{fig:Dimer_Prop_V}, Figure \ref{fig:Dimer_Prop_VS}, Figure \ref{fig:Dimer_Prop_ST}.


\section*{Discussion}

Adjoining a negative feedback strategy to the antithetic integral controller was shown to reduce the stationary variance for the controlled species, an effect that was expected from previous studies and predicted by the obtained theoretical results. The structure of the negative feedback strategy was notably emphasized to have important consequences on the magnitude of the variance reduction. Indeed, the ON/OFF controller can be used to dramatically reduce the variance while still preserving the ergodicity of the closed-loop network. This can be explained mainly because the proportional effective gain $\beta$ is very sensitive to changes in the feedback strength $K_p$ and can reach reasonably large values (still smaller than $K_p$); see Fig. \ref{fig:Gene_Prop_Beta} in the appendix. The preservation of the ergodicity property for the closed-loop network comes from the fact that $\E_\pi[K_p\max\{0,\mu-\theta X_\ell\}]$ remains smaller than the value of the nominal stationary control input (the constant input for which the stationary mean of the controlled species equals the desired set-point) for a wide range of values for $K_p$. Regarding the mean dynamics, this feedback leads to a decrease of the settling-time but also also leads to abrupt transient dynamics for large values of $K_p$ because of the presence of a stable zero in the mean closed-loop dynamics that is inversely proportional to $\beta$ (which is very sensitive to changes in $K_p$ in this case and which can reach high values). Unfortunately, this controller cannot be implemented in-vivo because it does not admit any reaction network implementation. However, it can still be implemented in-silico for the stochastic single-cell control for the control of cell populations using, for instance, targeted optogenetics; see e.g. \cite{Rullan:17}. On the other hand, the Hill feedback, while being practically implementable, has a much less dramatic impact on the stationary variance and on the mean dynamics. The first reason is that the effective proportional gain $\beta$ is less sensitive with respect to changes in $K_p$ and remains very small even when $K_p$ is large; see Fig. \ref{fig:Gene_Hill_Beta}. The absence of zero does not lead to any abrupt transient dynamics even for large values for $K_p$ but this may also be due to the fact that $\beta$ always remains small as opposed to as in the ON/OFF proportional feedback case. A serious issue with this feedback is that ergodicity can be easily lost since $\E_\pi[K_p/(1+X_\ell)]$ becomes very quickly larger than the value of the nominal control input as we increase $K_p$. The properties of both feedback strategies are summarized in Table \ref{tab:statsvsprop}.

To prove the main theoretical results, a tailored closure method had to be developed to deal with the bimolecular comparison reaction. A similar one has also been suggested in \cite{Olsman:17} for exactly the same purpose. These methods rely on the assumption that the molecular count of the controller species $\Z{2}$ is, most of the time, equal to 0, a property that is ensured by assuming  that $k/\eta\ll 1$. This allowed for the simplification and the closure of the moment equations. The theory was only developed for unimolecular networks because of the problem solvability. However, the extension of those theoretical results to more general reaction networks, such as bimolecular networks, is a difficult task mainly because of the moment closure problem that is now also present at the level of the species of the controlled network. In this regard, this extension is, at the moment, only possible using existing  moment closure methods (see e.g. \cite{Hespanha:08b,Milner:11,Smadbeck:13}) which are known to be potentially very inaccurate and would then compromise the validity of the obtained approximation. We believe that obtaining accurate and general theoretical approximations for the stationary variance for bimolecular networks is currently out of reach. It is also unclear whether the obtained qualitative and quantitative results still hold when the assumption $k/\eta\ll1$ on the controller parameters is not met.

Interestingly, the results obtained in the current paper provide some interesting insights on an unexpected connection between deterministic PI control and its stochastic analogues. In particular, it is possible to observe that the destabilizing effect of deterministic integral control is analogous to the variance increase due to the use of the stochastic antithetic integral controller. In a similar way, the stabilizing property of deterministic proportional controllers is the deterministic analogue of the property of variance decrease of the stochastic proportional controller; see Table \ref{tab:detvsstoc}).

The controller considered in this paper is clearly analogous to PI controllers. A usual complemental element is the so-called derivative action (or a filtered version of it) in order to add an anticipatory effect to the controller and prevent high overshoot; see \cite{Astrom:95}. So far, filtered versions of the derivative action have been proposed in a deterministic setting. Notably, the incoherent feedforward loop locally behaves like a filtered derivative action. More recently, a reaction network approximating a filtered derivative action was proposed in \cite{Halter:17} in the deterministic setting. It is unclear at the moment whether a stochastic version for the derivative action can be found but it is highly possible that such a stochastic derivative action can be implemented in terms of elementary reactions.

The negative feedback strategy considered here is an ideal/simplified one. Indeed, it was assumed in this paper that the controlled species was directly involved in the negative feedback. However, it is very likely that, the controlled species may not be directly usable in the feedback, that intermediary species may be involved (e.g. a gene expression network is involved in the feedback) or that the feedback is in terms of a species upstream the controlled species (for instance feedback uses a protein while the controlled species is the corresponding homodimer). The theory may be adapted to deal with such cases as long as the controlled network is unimolecular. It is however expected that the same qualitative behavior will be observed. The reason for that is that in unimolecular networks, species cooperate in the sense that they act positively on each other. Hence, decreasing the variance of one species will also decrease the variance of all the species that are created from it. For instance, in a gene expression network, if the mRNA variance is decreased, the protein variance will decrease as well, and vice-versa.

Finally, the implementation of such negative feedback loops is an important, yet elusive, task. It is unclear at the moment how in-vivo experiments could be conducted. Preliminary experimental results to validate the theoretical/computational ones could be obtained using optogenetics and single-cell control for population control. In-vivo experiments will certainly require a lot more effort.

\begin{table}
  \centering
  \caption{Effects of the different feedback strategies on the mean dynamics and the stationary variance.}\label{tab:statsvsprop}
  \begin{tabular}{|c||c|c|}
  \hline
                                    & ON/OFF Proportional Feedback & Hill Feedback\\
                                    \hline
                                    \hline
Ergodicity                 &  robust (+) & fragile (-)\\
\hline
$\beta$                  & very sensitive (+)  & poorly sensitive (-) \\
                                & wide range (+)  & small range \\
                                \hline
Mean Dynamics  &  reduce settling-time (+)   & reduce settling-time (+)\\
                                &  zero dynamics (-)   & no zero dynamics (+)\\
\hline
Stationary variance  & dramatic reduction (++)  &   slight reduction (+)\\
\hline
  \end{tabular}
\end{table}

\begin{table}
  \centering
  \caption{The effects of the proportional and integral actions on the dynamics of a system in both the deterministic and stochastic setting.}\label{tab:detvsstoc}
  \begin{tabular}{|c||c|c|}
  \hline
                                    & Integral action & Proportional action\\
                                    \hline
                                    \hline
    Deterministic  & regulation (+) & no regulation (-)\\
    Setting                            & destabilizing (-) & stabilizing (+)\\
                                \hline
    Stochastic         & regulation (+) & no regulation (-)\\
    Setting                & increases variance (-) & decreases variance (+)\\
                                \hline
  \end{tabular}
\end{table}

\appendix

\section*{Supplementary figures for the gene expression network}\label{secSI:gene_fig}

\begin{figure}[H]
  \centering
  \includegraphics[width=0.7\textwidth]{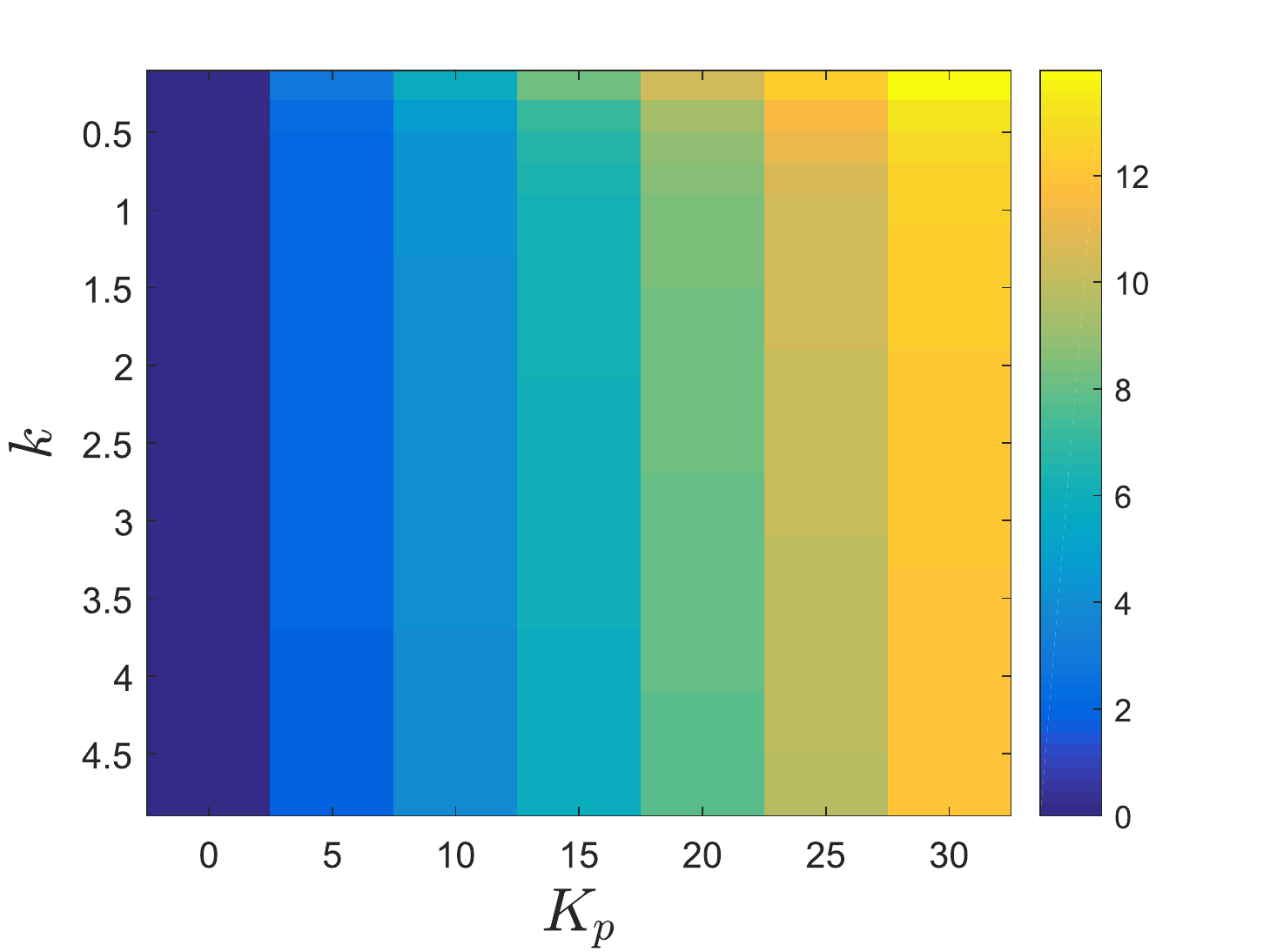}
  \caption{Evolution of $\beta$ as a function of the gains $k$ and $K_p$ in a gene network controlled with an antithetic integral controller combined with an ON/OFF proportional feedback.}\label{fig:Gene_Prop_Beta}
\end{figure}

\begin{figure}[H]
  \centering
  \includegraphics[width=0.7\textwidth]{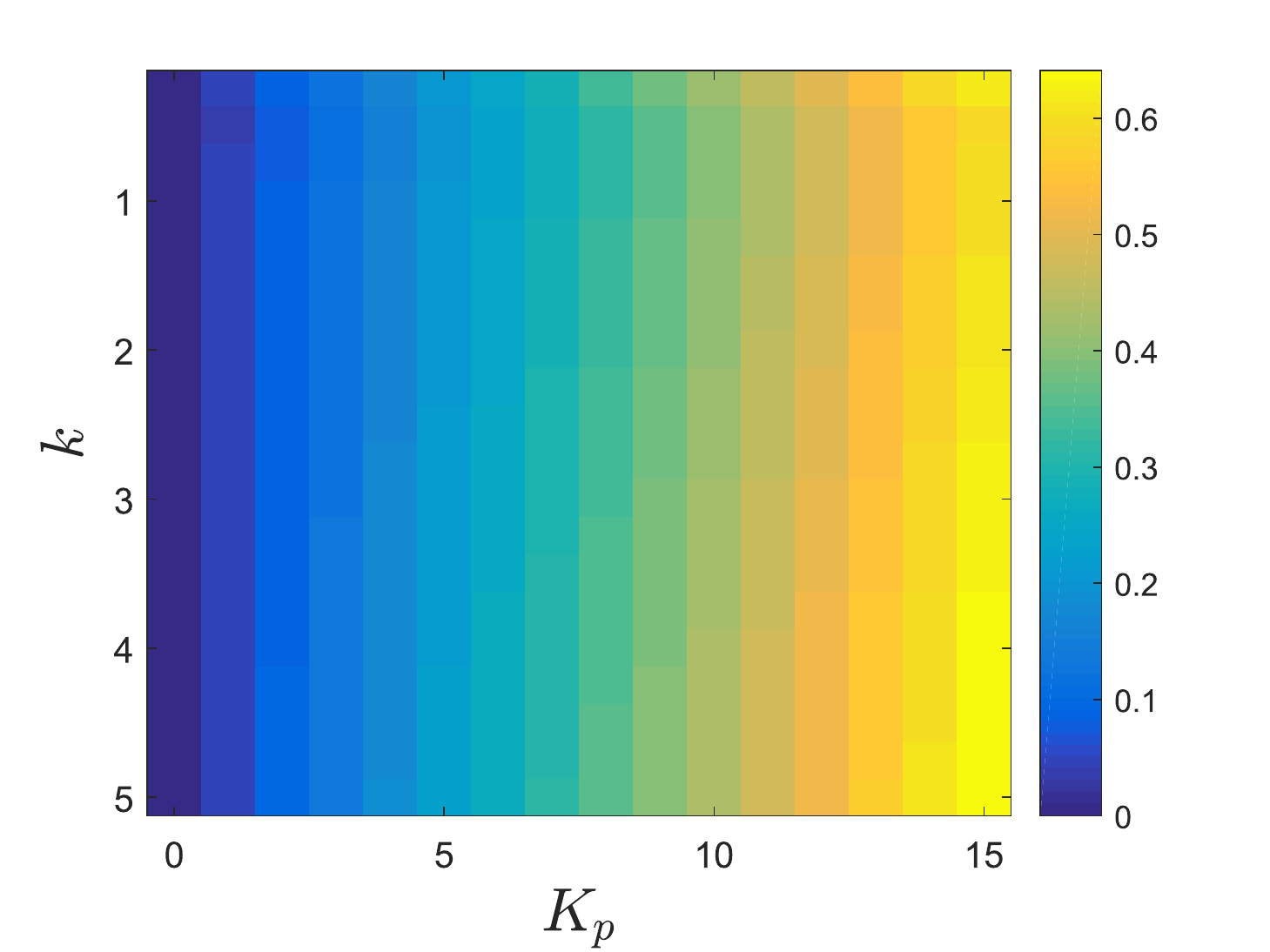}
  \caption{Evolution of $\beta$ as a function of the gains $k$ and $K_p$ in a gene network controlled with an antithetic integral controller combined with a Hill feedback.}\label{fig:Gene_Hill_Beta}
\end{figure}

\section*{Supplementary figures for the gene expression network with protein maturation}\label{secSI:genem_fig}

\begin{figure}[H]
  \centering
  \includegraphics[width=0.7\textwidth]{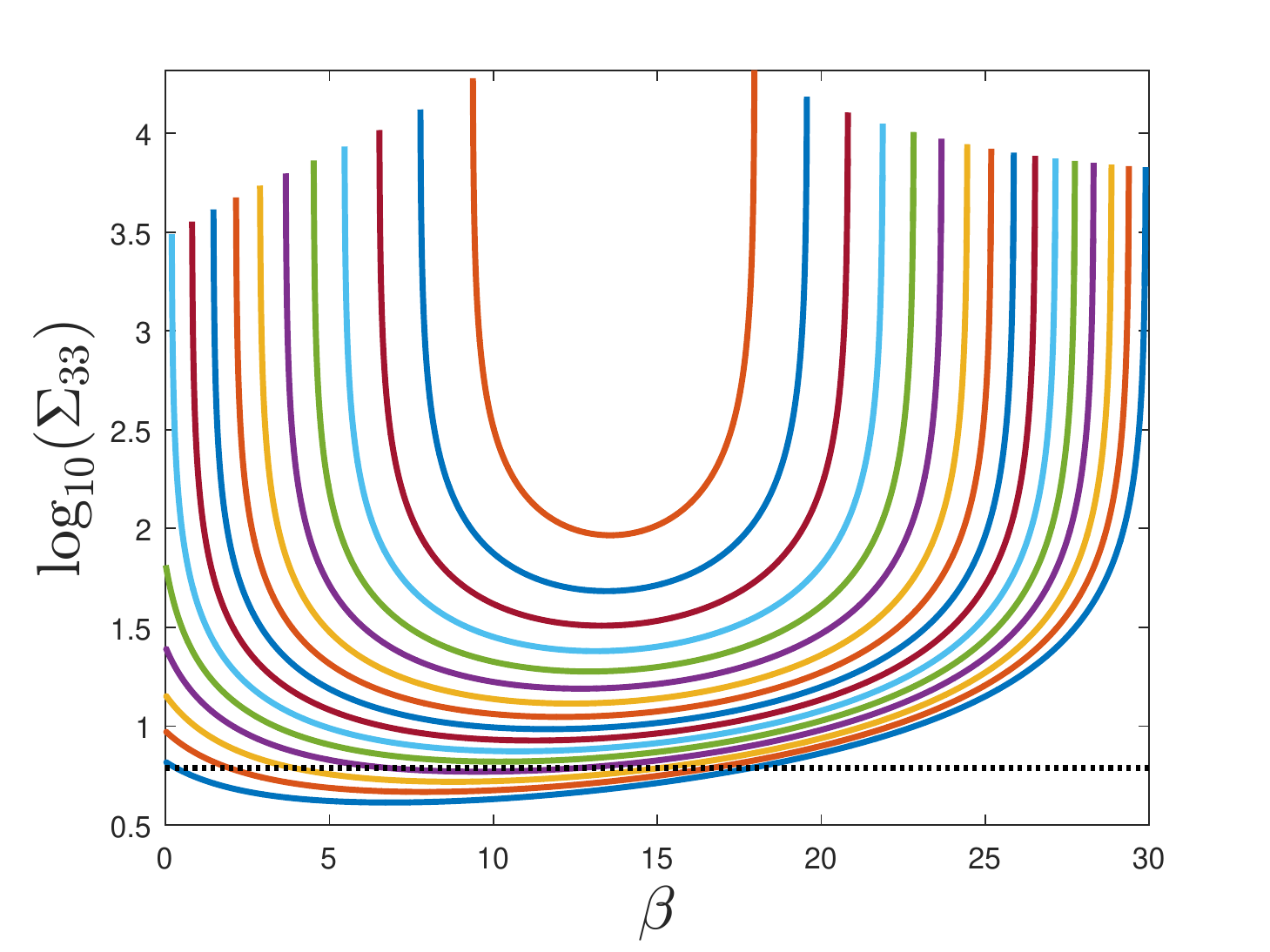}
  \caption{Evolution of the logarithm of the predicted stationary variance for the mature protein copy number as a function of  the gain $\beta$. We can observe a nonmonotonic behavior due to the fact that the stability region for $R$ in terms of the parameters $k$ and $\beta$ is bounded. The black dashed line is the logarithm of the open-loop variance.}\label{fig:Mat_NM}
\end{figure}

\begin{figure}[H]
  \centering
  \includegraphics[width=0.7\textwidth]{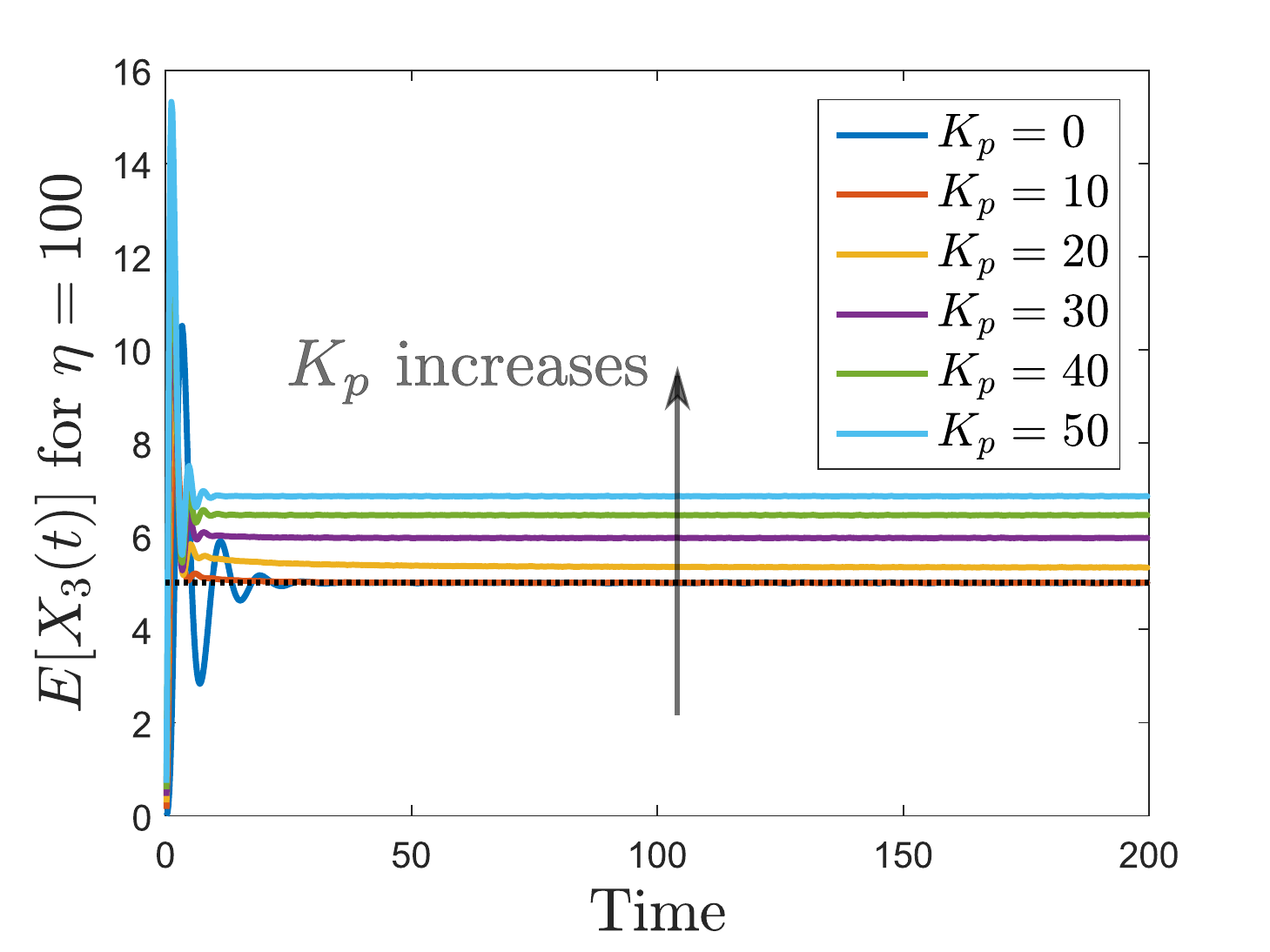}
 \caption{Mean trajectories for the mature protein copy number when the gene expression network with protein maturation is controlled with the antithetic integral controller \eqref{eq:AIC} with $k=3$ and an ON/OFF proportional controller. The set-point value is indicated as a black dotted line.}\label{fig:Mat_Prop_E_NM}
\end{figure}

\begin{figure}[H]
  \centering
  \includegraphics[width=0.7\textwidth]{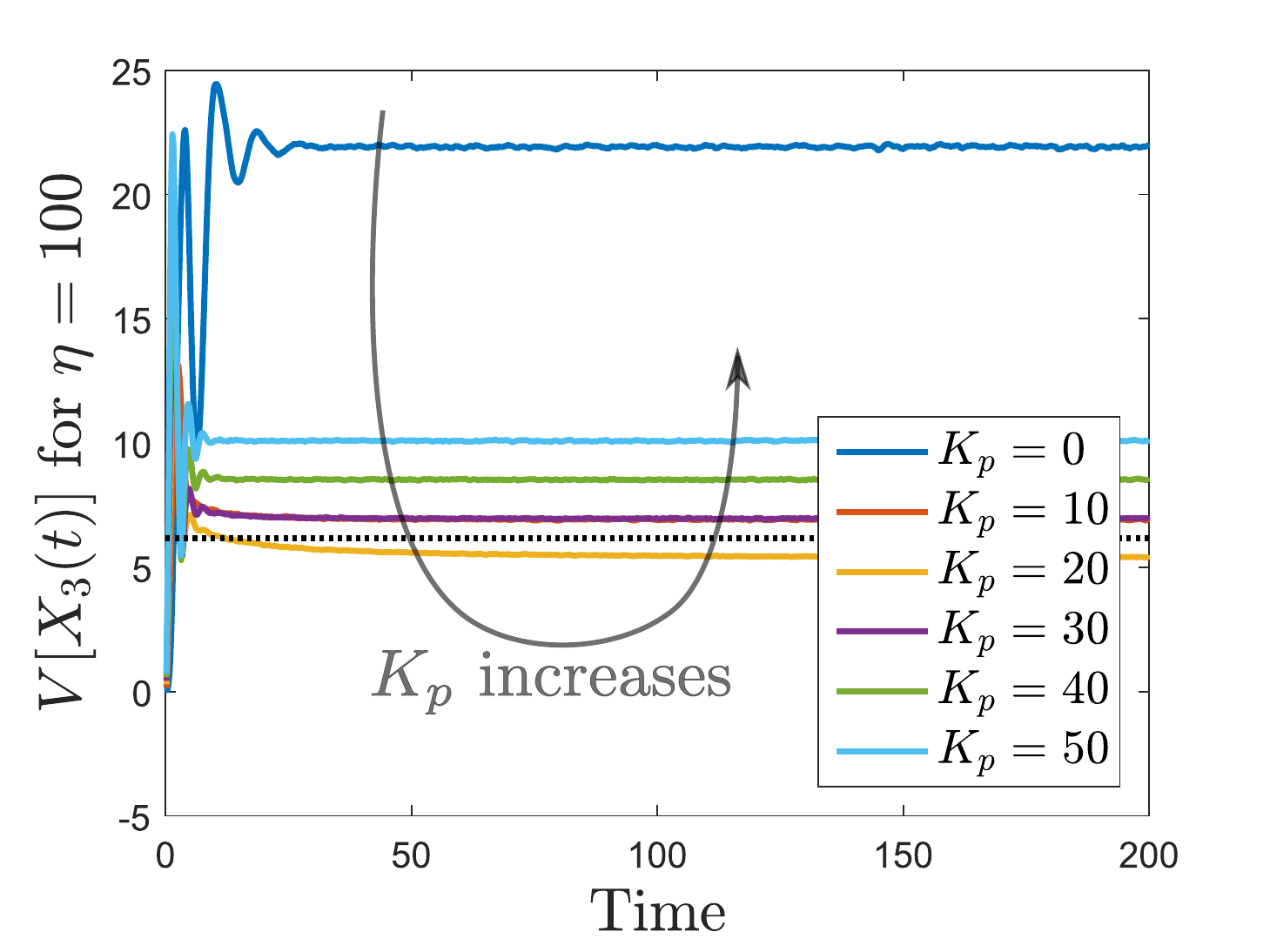}
 \caption{Variance trajectories for the mature protein copy number when the gene expression network with protein maturation is controlled with the antithetic integral controller \eqref{eq:AIC} with $k=3$ and an ON/OFF proportional  controller.  The stationary constitutive variance is depicted in black dotted line.}\label{fig:Mat_Prop_V_NM}
\end{figure}

\begin{figure}[H]
  \centering
  \includegraphics[width=0.7\textwidth]{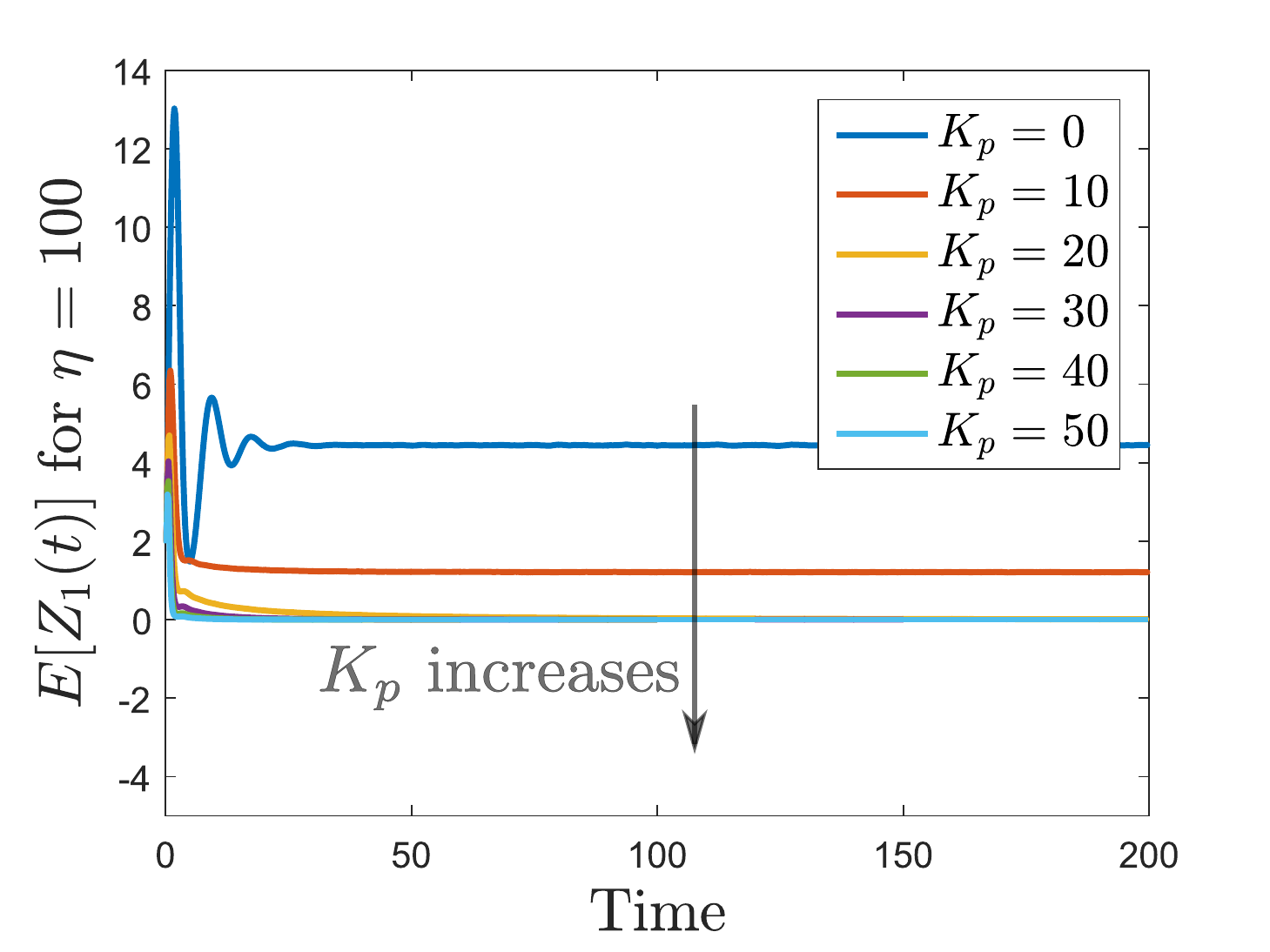}
 \caption{Mean trajectories for the actuating species copy number when the gene expression network with protein maturation is controlled with the antithetic integral controller \eqref{eq:AIC} with $k=3$ and an ON/OFF proportional controller. }\label{fig:Mat_Prop_Z1_NM}
\end{figure}

\begin{figure}[H]
  \centering
  \includegraphics[width=0.7\textwidth]{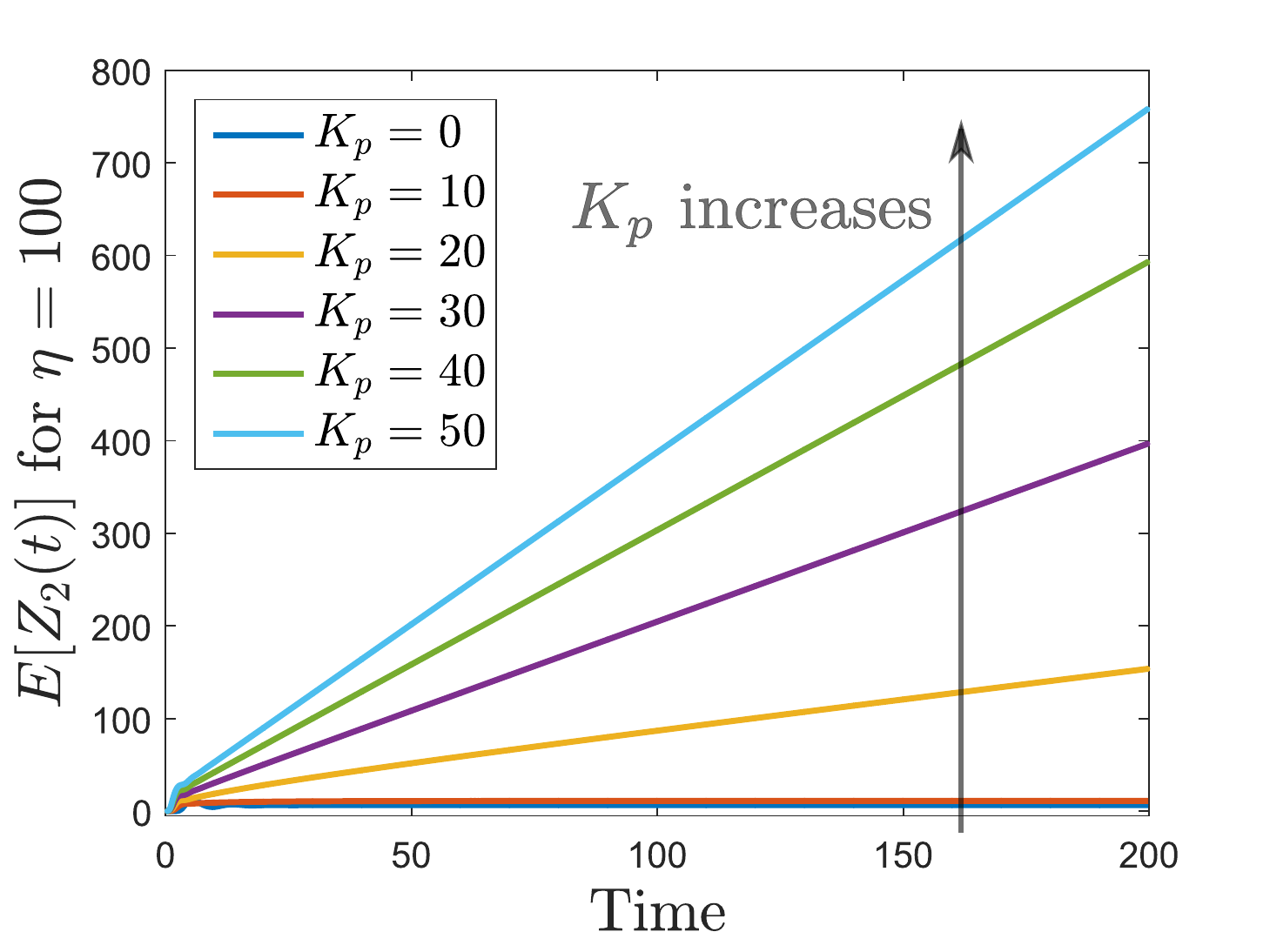}
 \caption{Mean trajectories for the sensing species copy number when the gene expression network with protein maturation is controlled with the antithetic integral controller \eqref{eq:AIC} with $k=3$ and an ON/OFF proportional  controller.}\label{fig:Mat_Prop_Z2_NM}
\end{figure}

\begin{figure}[H]
  \centering
  \includegraphics[width=0.7\textwidth]{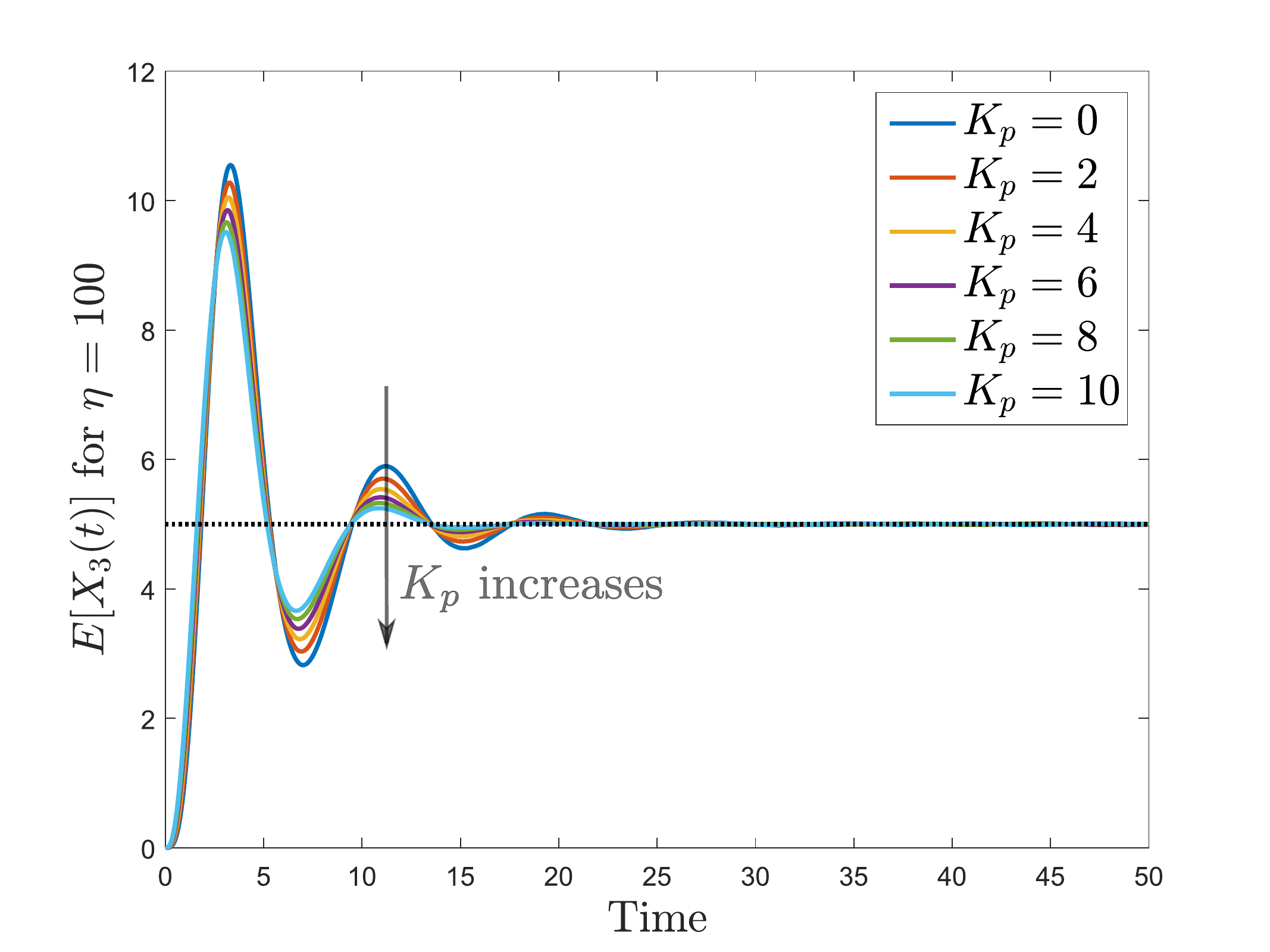}
 \caption{Mean trajectories for the mature protein copy number when the gene expression network with protein maturation is controlled with the antithetic integral controller \eqref{eq:AIC} with $k=3$ and a Hill negative feedback controller. The set-point value is indicated as a black dotted line.}\label{fig:Mat_Hill_E}
\end{figure}

\begin{figure}[H]
  \centering
  \includegraphics[width=0.7\textwidth]{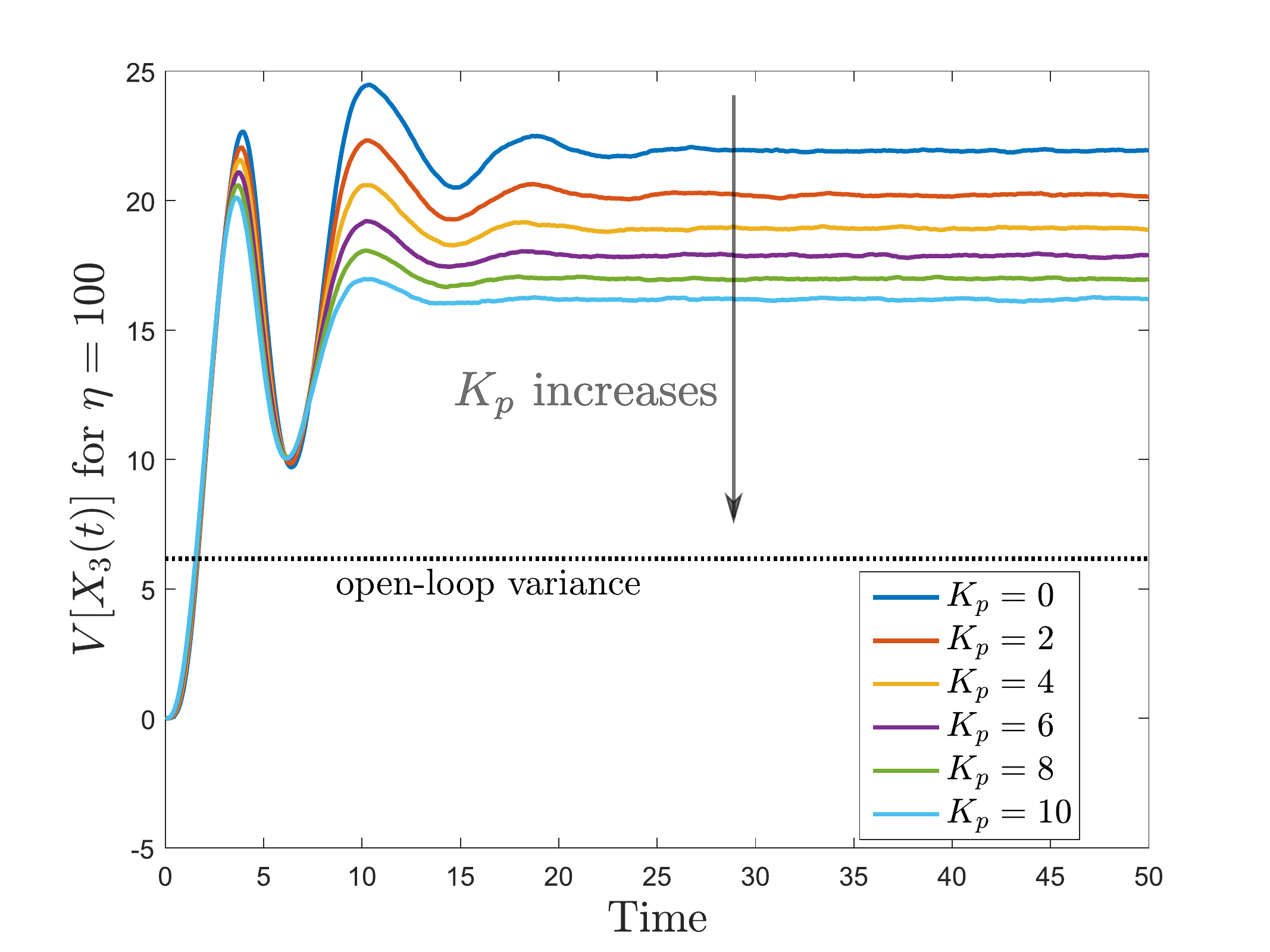}
 \caption{Variance trajectories for the mature protein copy number when the gene expression network with protein maturation is controlled with the antithetic integral controller \eqref{eq:AIC} with $k=3$ and a Hill negative feedback controller.  The stationary constitutive variance is depicted in black dotted line.}\label{fig:Mat_Hill_V}
\end{figure}

\begin{figure}[H]
  \centering
  \includegraphics[width=0.7\textwidth]{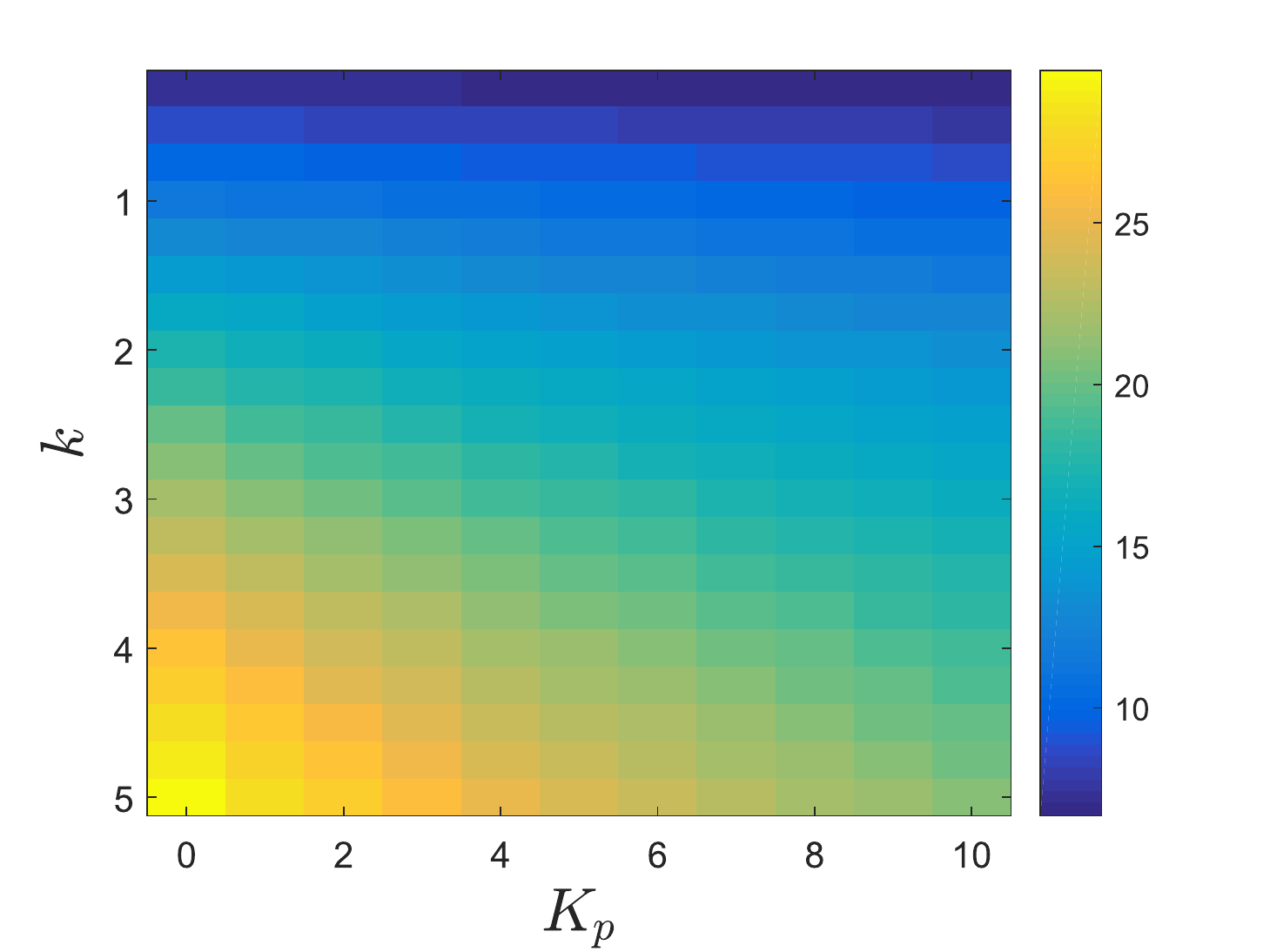}
  \caption{Stationary variance for the mature protein copy number when the gene expression network with protein maturation is controlled with the antithetic integral controller \eqref{eq:AIC} and a Hill feedback controller.}\label{fig:Mat_Hill_VS}
\end{figure}

\begin{figure}[H]
  \centering
  \includegraphics[width=0.7\textwidth]{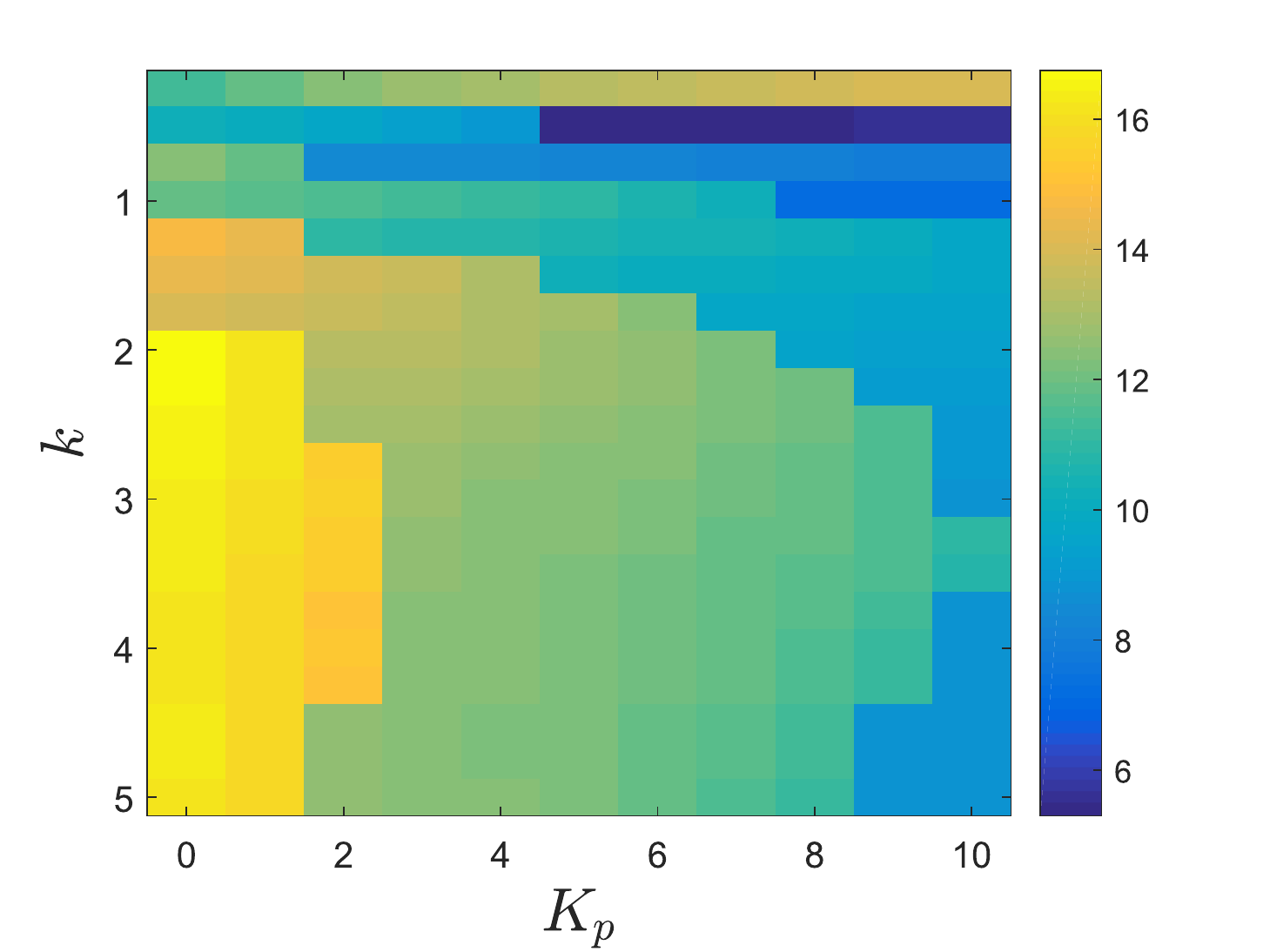}
   \caption{Settling-time for the mean trajectories for the mature protein copy number when the gene expression network with protein maturation is controlled with the antithetic integral controller \eqref{eq:AIC} and a Hill negative feedback controller.}\label{fig:Mat_Hill_ST}
\end{figure}

\section*{Supplementary figures for the gene expression network with protein dimerization}\label{secSI:gened_fig}

\begin{figure}[H]
  \centering
  \includegraphics[width=0.7\textwidth]{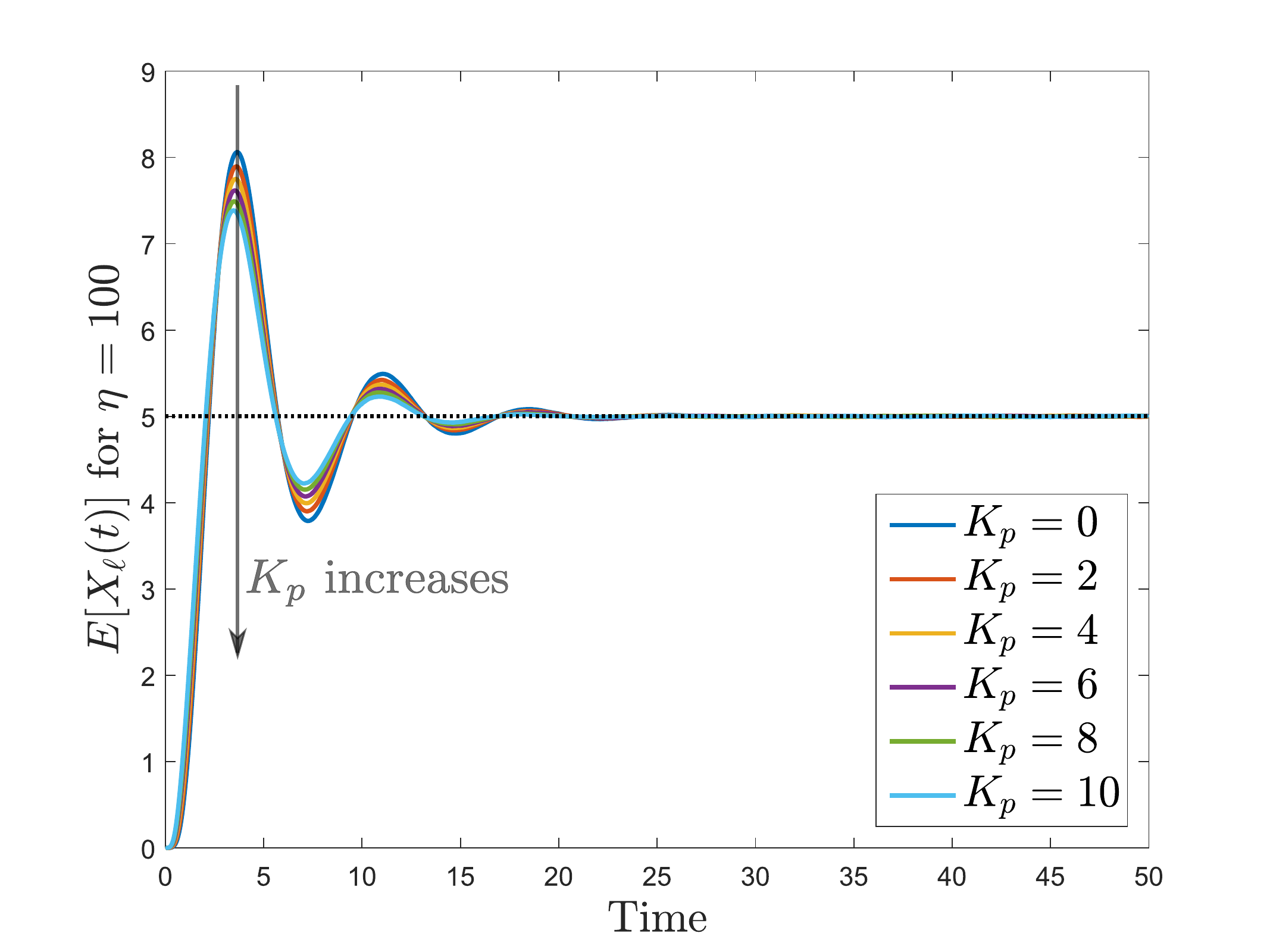}
    \caption{Mean trajectories for the homodimer copy number when the gene expression network with protein dimerization is controlled with the antithetic integral controller \eqref{eq:AIC} with $k=3$ and a Hill negative feedback. The set-point value is indicated as a black dotted line.}\label{fig:Dimer_Hill_E}
\end{figure}

\begin{figure}[H]
  \centering
  \includegraphics[width=0.7\textwidth]{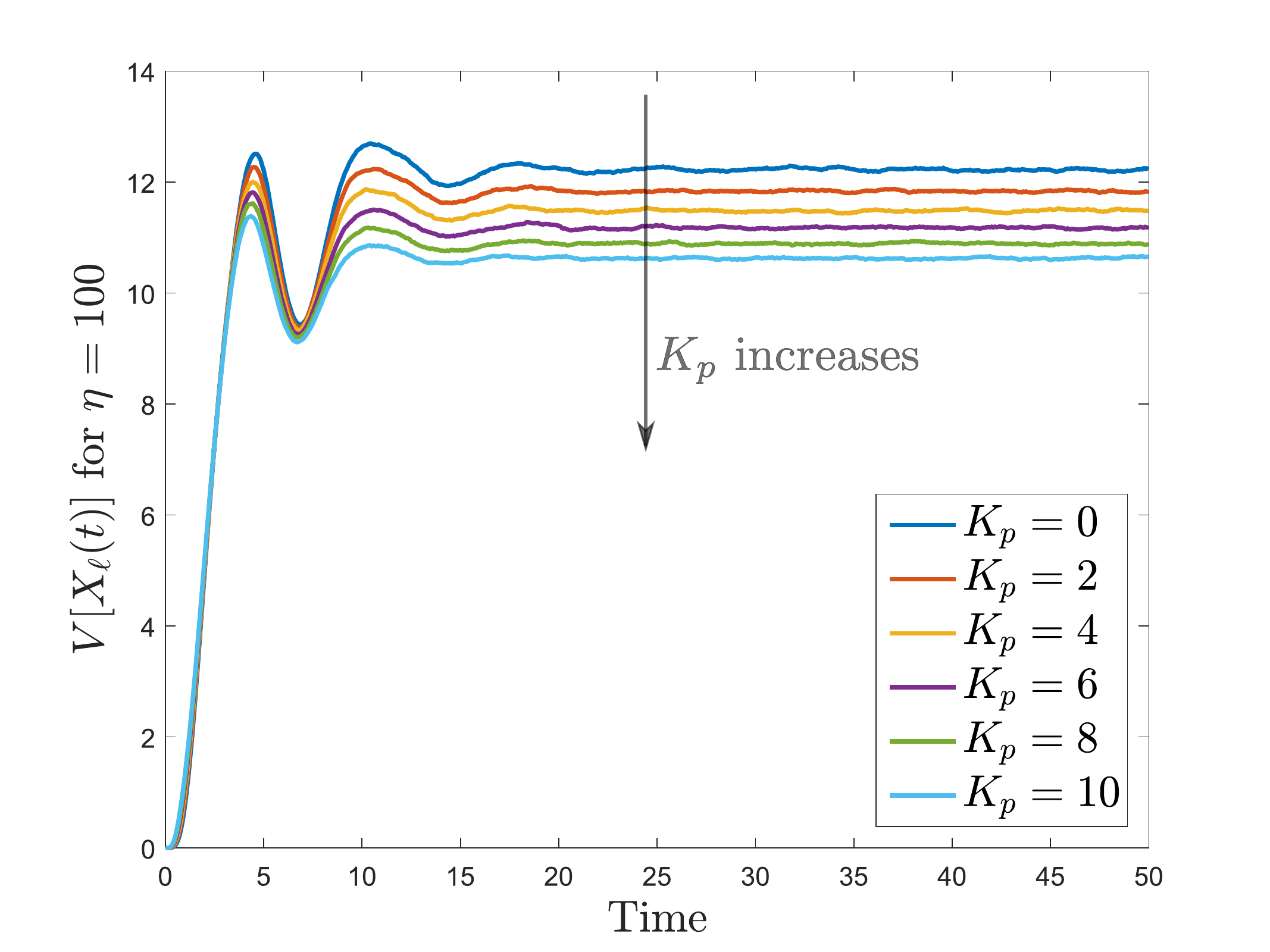}
  \caption{Variance trajectories for the homodimer copy number when the gene expression network with protein dimerization  is controlled with the antithetic integral controller \eqref{eq:AIC} with $k=3$ and a Hill negative feedback.  The stationary constitutive variance is depicted in black dotted line.}\label{fig:Dimer_Hill_V}
\end{figure}

\begin{figure}[H]
  \centering
  \includegraphics[width=0.7\textwidth]{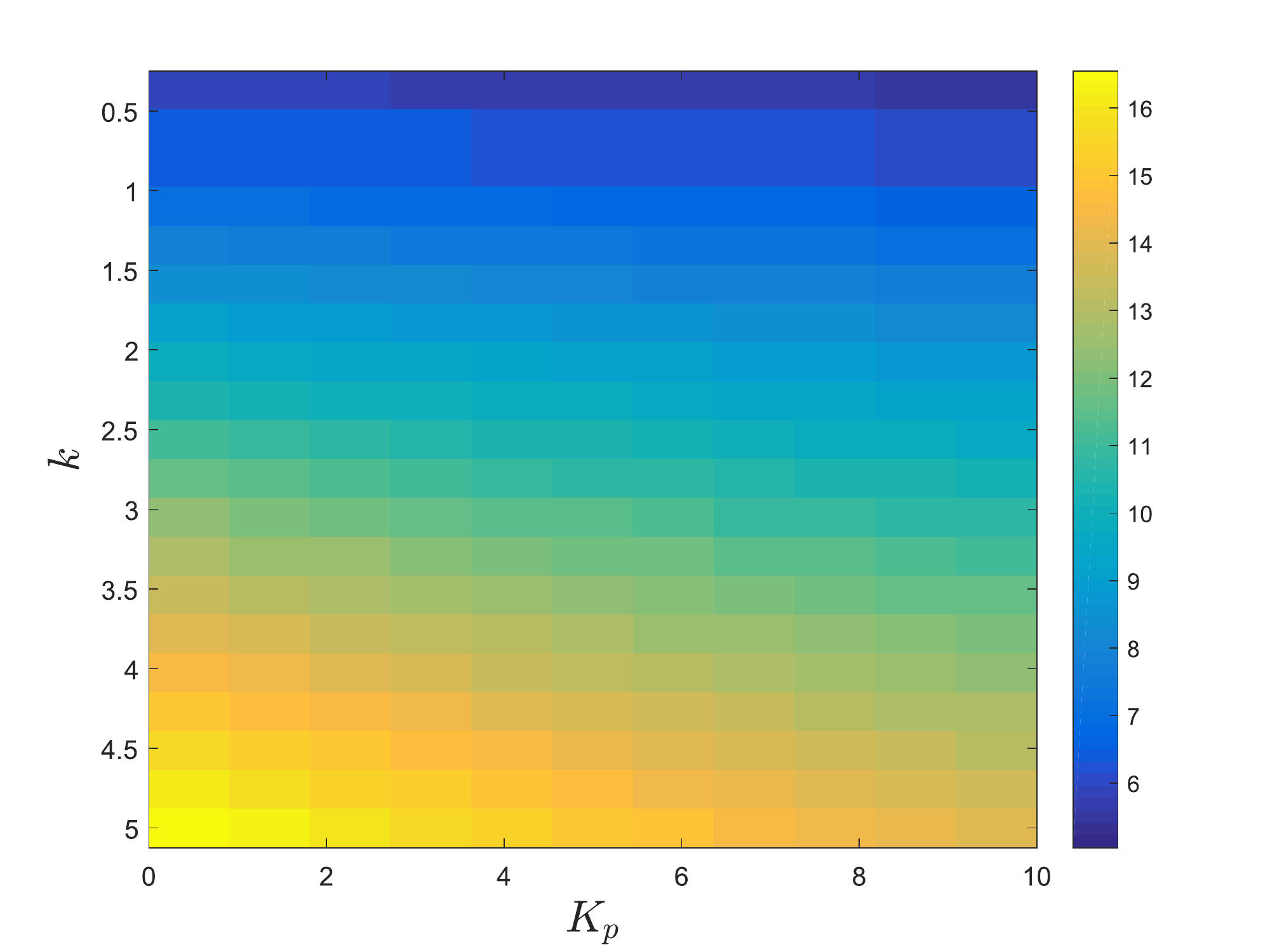}
   \caption{Stationary variance for the homodimer copy number when the gene expression network with protein dimerization  is controlled with the antithetic integral controller \eqref{eq:AIC} and aa Hill negative feedback. }\label{fig:Dimer_Hill_VS}
\end{figure}

\begin{figure}[H]
  \centering
  \includegraphics[width=0.7\textwidth]{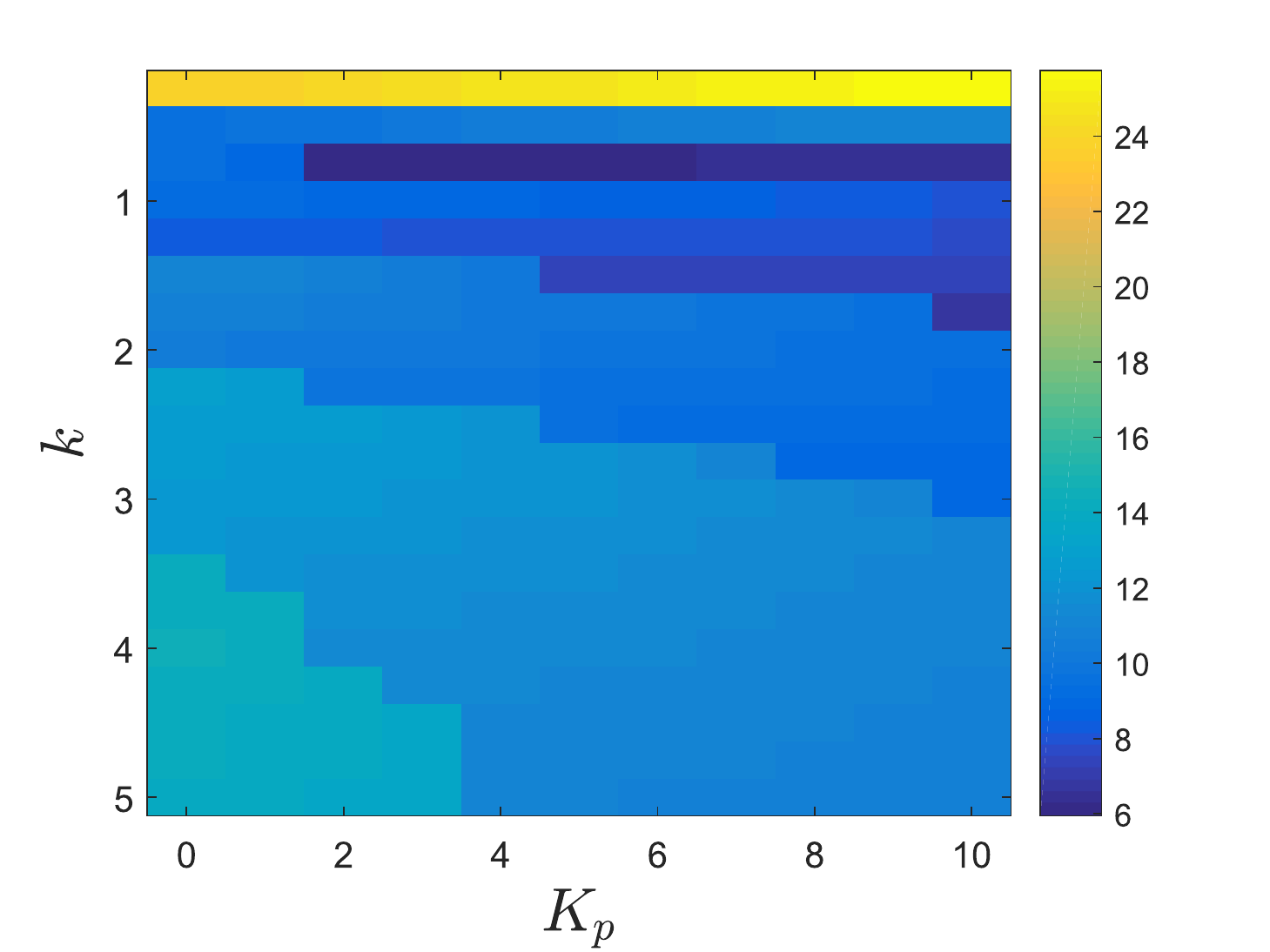}
  \caption{Settling-time for the mean trajectories for the homodimer copy number when the gene expression network with protein dimerization  is controlled with the antithetic integral controller \eqref{eq:AIC} and a Hill negative feedback. }\label{fig:Dimer_Hill_ST}
\end{figure}

\bibliographystyle{plain}

\end{document}